\documentclass[review]{elsarticle}                                     
\usepackage{lineno,hyperref}
\modulolinenumbers[5]

\usepackage{graphicx}          
\usepackage{amsmath,graphicx} 
\usepackage{amssymb}
\usepackage[dvips]{epsfig}
\usepackage{subfigure}
\usepackage{natbib}
\usepackage{multirow}
\usepackage{mathrsfs}
\usepackage{enumerate}
\usepackage{color}
\usepackage{mathtools}

\DeclareMathOperator*{\argmin}{arg\,min}

\begin{document}

\begin{frontmatter}

\title{Receding-horizon Stochastic Model Predictive Control with Hard Input Constraints and Joint State Chance Constraints}

\author[MIT,Berkeley]{Joel A. Paulson}
\author[Berkeley]{Edward A. Buehler}
\author[MIT]{Richard D. Braatz}
\author[Berkeley]{Ali Mesbah}
\ead{mesbah@berkeley.edu}

\address[MIT]{Department of Chemical Engineering, Massachusetts Institute of Technology, Cambridge, MA 02139, USA}
\address[Berkeley]{Department of Chemical and Biomolecular Engineering, University of California-Berkeley, Berkeley, CA 94720, USA}
          
\begin{abstract}                          

This article considers the stochastic optimal control of discrete-time linear systems subject to (possibly) unbounded stochastic disturbances, hard constraints on the manipulated variables, and joint chance constraints on the states. A tractable convex second-order cone program (SOCP) is derived for calculating the receding-horizon control law at each time step. Feedback is incorporated during prediction by parametrizing the control law as an affine function of the disturbances. Hard input constraints are guaranteed by saturating the disturbances that appear in the control law parametrization. The joint state chance constraints are conservatively approximated as a collection of individual chance constraints that are subsequently relaxed via the Cantelli-Chebyshev inequality. Feasibility of the SOCP is guaranteed by softening the approximated chance constraints using the exact penalty function method. Closed-loop stability in a stochastic sense is established by establishing that the states satisfy a geometric drift condition outside of a compact set such that their variance is bounded at all times. The SMPC approach is demonstrated using a continuous acetone-butanol-ethanol fermentation process, which is used for production of high-value-added drop-in biofuels. 
\end{abstract}

\begin{keyword}
Stochastic systems \sep Predictive control \sep Constrained control \sep Continuous acetone-butanol-ethanol fermentation  
\end{keyword}

\end{frontmatter}


\section{Introduction}
\label{sec:intro}

Robust model predictive control (MPC) has been the subject of intensive study over the past two decades (e.g., see \cite{bemporad99,lim09,mayne14} and the references therein). Deterministic approaches to robust MPC consider uncertainties of bounded nature, and typically rely on worst-case optimal control formulations in which system constraints are satisfied for all possible uncertainty realizations. Designing a control policy with respect to all uncertainty realizations may, however, be overly conservative since worst-case uncertainty realizations often have a low likelihood of occurrence. In addition, it is often impractical to precisely determine uncertainty bounds for use in deterministic robust MPC approaches.

Stochastic MPC (SMPC) has recently gained increasing interest to reduce the conservatism of deterministic robust MPC approaches by explicitly incorporating probabilistic descriptions of uncertainties into the optimal control problem. Also, SMPC can take advantage of the fact that, in certain applications, system constraints need not be satisfied for all uncertainty realizations. In such systems, hard system constraints can be replaced with chance constraints (or, alternatively, expectation-type constraints) to allow for an admissible level of constraint violation. This formulation enables systematically trading off robustness to uncertainties (in terms of constraint satisfaction) with closed-system performance.                    

Broadly speaking, SMPC approaches for stochastic linear systems can be categorized into three classes (see \cite{mesbah15} for a recent review on SMPC). The first class consists of \textsl{stochastic tube approaches} \cite{cannon09,cannon11,kouv15} that use stochastic tubes with fixed or variable cross sections to replace chance constraints with linear constraints on the nominal state predictions and to construct terminal sets for guaranteeing recursive feasibility. These approaches use a prestabilizing feedback controller to ensure closed-loop stability. However, since the prestabilizing state feedback controller is determined offline (only open-loop control actions are used as decision variables in the online optimization), stochastic tube approaches cannot handle hard input constraints. Less suboptimal\footnote{Suboptimality is considered with respect to dynamic programming (DP), which provides the most general framework for control of uncertain systems \cite{lee04}. DP considers arbitrary control laws, which commonly makes it impractical for real control applications due to \textsl{curse of dimensionality}.}SMPC formulations entail online computation of both the feedback control gains and the open-loop control actions. Such formulations, however, lead to a nonconvex optimization \cite{farina15}. 

The second class of SMPC approaches use \textsl{affine parameterizations of feedback control laws} to obtain a tractable convex SMPC formulation \cite{oldew08,hokayem09,hokayem12,paulson15}. A key challenge in these approaches is to handle hard input constraints in the presence of unbounded stochastic uncertainties (e.g., Gaussian noise), as unbounded uncertainties almost surely lead to excursions of states from any bounded set. To address this challenge, \cite{hokayem09} introduced \textsl{saturation functions} into affine feedback control policies. Saturation functions allow for rendering the feedback control policies nonlinear to directly handle hard input constraints, without the need to relax the hard input constraints to soft chance constraints. The third class of SMPC approaches entail the so-called \textsl{sample-based approaches} (aka scenario approaches), which essentially characterize the stochastic system dynamics using a finite set of random realizations of uncertainties that are used to solve the optimal control problem in one shot \cite{batina04,blackmore10,calaf13}. Even though sample-based SMPC approaches commonly do not rely on any convexity requirements, establishing recursive feasibility and closed-loop stability for these approaches is generally challenging, particularly for the case of unbounded uncertainties.

This paper considers the problem of receding-horizon control of stochastic linear systems with (possibly) unbounded disturbances that have arbitrary distributions. A SMPC approach is presented that can handle joint state constraints and hard input constraints in the presence of the unbounded disturbances. The primary challenges in solving the SMPC problem include general intractability of the problem for arbitrary feedback control laws, unbounded nature of disturbances, and nonconvexity of chance constraints. These challenges are addressed by using a saturated affine disturbance feedback control policy and by approximating the joint state chance constraints as a set of individual chance constraints that can then be relaxed via the Cantelli-Chebyshev inequality. A tractable convex formulation is derived for the SMPC problem in terms of a second order cone program with guaranteed feasibility and stability. 

The problem setup considered in this paper is similar to that in \cite{hokayem12,farina15}. Hard input constraints are relaxed to chance constraints in \cite{farina15} to cope with the unbounded nature of disturbances, whereas the SMPC program presented in this work allows for directly handling the hard input constraints in an unbounded disturbance setting. In addition, the SMPC approach in \cite{farina15} results in a nonconvex program, which can become impractical for large-scale systems. On the other hand, the main contribution of this paper with respect to \cite{hokayem12} is the ability to handle joint chance constraints while retaining the hard input constraints. Furthermore, unlike \cite{hokayem12}, the presented SMPC approach guarantees stability \textsl{by design} without using a stability constraint that can possibly limit the domain of attraction of the controller. In this work, the SMPC approach is demonstrated using an experimentally validated continuous acetone-butanol-ethanol (ABE) fermentation case study \cite{haus11}, which consists of $12$ state variables. The ABE fermentation process is used for production of high-value-added drop-in biofuels from lignocellulosic biomass.                          

The organization of the paper is as follows. The problem setup is outlined in Section~\ref{sec:problem_statement}. Section~\ref{sec:methods} discusses the methods used for deriving the tractable SMPC program presented in Section~\ref{sec:SMPC}, which is followed by the continuous ABE fermentation case study in Section \ref{sec:case_study} and conclusions in Section \ref{sec:conclusions}. For completeness, the proofs of all theoretical results presented in this paper are included in the Appendix A.

\textbf{Notation.} Hereafter, $\mathbb{N} = \{1,2,\ldots\}$ is the set of natural numbers. $\mathbb{N}_0 \triangleq \mathbb{N} \cup \{0\}$. $\mathbb{Z}_{[a,b]}:=\{a,a+1,\ldots,b\}$ is the set of integers from $a$ to $b$. $\mathbb{R}_{+}$ and $\mathbb{R}_{++}$ are the set of nonnegative and positive real numbers, respectively. $\mathbb{S}^n \subset \mathbb{R}^{n\times n}$ is the set of real symmetric matrices. $\mathbb{S}^n_+$ and $\mathbb{S}^n_{++}$ are the set of positive semidefinite and definite matrices, respectively. $\emptyset \triangleq \{ \}$ is the empty set. $I_a$ is the $a\times a$ identity matrix. tr$(\cdot)$ is the trace of a square matrix. $\| \cdot \|_p$ and $\| \cdot \|$ are the standard $p$-norm and the Euclidean norm, respectively. $\| x \|^2_A \triangleq x^\top A x$ is the weighted 2-norm. $\lambda_{\text{min}}(\cdot)$ and $\lambda_{\text{max}}(\cdot)$ are the minimum and maximum eigenvalues of a matrix, respectively. $\text{diag}(\cdot)$ is a block diagonal matrix. $\otimes$ denotes the Kronecker product. $\mathbf{E}[\cdot]$ or $\bar{(\cdot)}$ denotes the expected value. $\mathbf{Var}[\cdot]$ or $\Sigma_{(\cdot)}$ denotes the covariance matrix. $\mathbf{Pr}[\cdot]$ denotes probability. $\mathbf{E}_x[\cdot]$, $\mathbf{Var}_x[\cdot]$, and $\mathbf{Pr}_x[\cdot]$ denote the conditional expected value, conditional variance, and conditional probability given information $x$, respectively. $\sigma_{a,b} \triangleq \mathbf{E}[(a-\bar{a})(b-\bar{b})^\top]$ is the cross-covariance between random vectors $a$ and $b$. $\mathbf{1}_A(\cdot)$ is the indicator function of the set $A$. $\mathcal{N}(\mu,\Sigma)$ is the Gaussian distribution with mean $\mu$ and covariance $\Sigma$. 

\section{Problem Statement}
\label{sec:problem_statement}

Consider a discrete-time, stochastic linear system
\begin{align} \label{eq:plantmodel}
x^+ = A x + B u + G w,
\end{align}
where $x \in \mathbb{R}^{n_x}$, $u \in \mathbb{R}^{n_u}$, and $w \in \mathbb{R}^{n_w}$ are the system states, inputs, and disturbances at the current time instant, respectively; $x^+$ denotes the system states at the next time instant; and $A$, $B$, and $G$ are the known system matrices. Note that system~\eqref{eq:plantmodel} is a Markov process. The following assumptions are made throughout the paper. 

\textbf{Assumption 1}
\begin{itemize}
	\item[(i)] \textit{The system matrix $A$ is Schur stable, i.e., all of its eigenvalues lie inside the unit circle in the complex plane.}
	\item[(ii)] \textit{The states $x$ can be observed exactly at all times.}
	\item[(iii)] \textit{The stochastic disturbances $w$ are assumed to be zero-mean independent and identically distributed (i.i.d.) random variables with known covariance matrix $\Sigma_w \in \mathbb{S}^{n_w}_{++}$ and arbitrary, but known probability density function (PDF) $p_w$. The (possibly) unbounded support of $p_w$ is denoted by $\Delta_w \subseteq \mathbb{R}^{n_w}$.}
	\item[(iv)] \textit{The system is subject to hard input constraints of the polytopic form}
\begin{align} \label{eq:inputcons}
u \in \mathbb{U} \triangleq \{ u \in \mathbb{R}^{n_u} \mid H_u u \leq k_u \},
\end{align}
\textit{where $H_u \in \mathbb{R}^{s \times n_u}$; $k_u \in \mathbb{R}^s$; and $s \in \mathbb{N}$ is the number of input constraints. Further, $\mathbb{U}$ is assumed to be a compact set that contains the origin in its interior. Note that \eqref{eq:inputcons} can represent both $1$- and $\infty$-norm sets.}
	\item[(v)] \textit{The states must satisfy joint chance constraints (JCCs) of the form}
\begin{align} \label{eq:statecc}
\mathbf{Pr}[x \in \mathbb{X}] \geq 1-\beta,
\end{align}
\textit{where $\beta \in [0,1)$ is the maximum allowed probability of constraint violation; $\beta = 0$ corresponds to the case where state constraints should hold for all disturbance realizations. The set},
\begin{align} \label{eq:stateconsset}
\mathbb{X}\triangleq \{ x \in \mathbb{R}^{n_x} \mid H_x x \leq k_x \}
\end{align}
\textit{is polytopic, where $H_x \in \mathbb{R}^{r \times n_x}$, $k_x \in \mathbb{R}^{r}$, and $r \in \mathbb{N}$ is the number of state constraints.} $\hfill \triangleleft$ 
\end{itemize}

This work considers the design of a SMPC approach for the stochastic system~\eqref{eq:plantmodel} such that it can handle hard input constraints~\eqref{eq:inputcons} and state chance constraints~\eqref{eq:statecc}. In addition, the control approach should guarantee stability of the closed-loop system in a stochastic sense. In the sequel, the stochastic optimal control problem pertaining to the aforementioned system setup is formulated. 

Let $N \in \mathbb{N}$ denote the prediction horizon of the control problem, and define $\mathbf{w}\triangleq [ w_0^\top,w_1^\top,\ldots,w_{N-1}^\top]^\top$ to be a disturbance sequence over the time interval $0$ to $N-1$. A general full-state \textsl{feedback control policy} $\boldsymbol\pi$ over the horizon $N$\footnote{For notational convenience, the control horizon is assumed to be equal to the prediction horizon.} is defined by
\begin{align}  \label{eq:input}
\boldsymbol\pi \triangleq \{ \pi_0,\pi_1(\cdot)\ldots,\pi_{N-1}(\cdot) \},
\end{align}
where $\pi_0 \in \mathbb{U}$ are the control inputs applied to the system~\eqref{eq:plantmodel} (i.e., $u=\pi_0$); and $\pi_i(\cdot) : \mathbb{R}^{n_x} \rightarrow \mathbb{U}, \, \forall \, i \in \mathbb{Z}_{[1,N-1]}$ are feedback control laws that are functions of future states. Let $\tilde{x}_i(x_0,\boldsymbol\pi,\mathbf{w})$ be state predictions for the system \eqref{eq:plantmodel} at time $i \in \mathbb{Z}_{[0,N]}$. $\tilde{x}_i(x_0,\boldsymbol\pi,\mathbf{w})$ is obtained by solving \eqref{eq:plantmodel} when the initial states are $x_0$ at time instant $0$, the control laws $\pi_j$ are applied at times $j \in \mathbb{Z}_{[0,i-1]}$, and the disturbance realizations are $\{w_0,\ldots,w_{i-1}\}$. The prediction model is written as 
\begin{align} \label{eq:statepred}
\tilde{x}_{i+1} = A\tilde{x}_i + B\pi_i(\tilde{x}_i) + Gw_i.
\end{align}
In the remainder of the paper, the explicit functional dependencies of $\{ \tilde{x}_i(x_0,\boldsymbol\pi,\mathbf{w}) \}_{i=0}^N$ on the initial states $x_0$, feedback control policy $\boldsymbol\pi$, and disturbance sequence $\mathbf{w}$ are dropped.

For the linear system \eqref{eq:plantmodel} with unbounded stochastic disturbances, the receding-horizon SMPC problem with hard input constraints and joint state chance constraints can now be stated as follows. 

\textbf{Problem 1 (Receding-horizon SMPC).} \textit{Given the current states $x$ observed from the system \eqref{eq:plantmodel}, receding-horizon SMPC entails solving the stochastic optimal control problem at each sampling time instant:}
\begin{align} \label{e_P1}
&\min_{\boldsymbol\pi} ~ V_N(x,\boldsymbol\pi) \\\notag
& \begin{array}{lll}
\text{s.t.:} & \tilde{x}_{i+1} = A\tilde{x}_i + B\pi_i(\tilde{x}_i) + Gw_i, &\forall i\in\mathbb{Z}_{[0,N-1]} \\[-1mm]
& \mathbf{Pr}[\pi_i(\tilde{x}_i) \in \mathbb{U}] = 1, &\forall  i\in\mathbb{Z}_{[0,N-1]} \\[-1mm]
& \mathbf{Pr}\big[ \tilde{x}_i \in \mathbb{X} \big] \geq 1-\beta, &\forall i\in\mathbb{Z}_{[1,N]} \\[-1mm]
& \tilde{x}_0 = x,~\mathbf{w} \sim p^N_w,
\end{array}
\end{align}
\textit{where the value function $V_N(x,\boldsymbol\pi)$ is defined by}
\begin{align} \label{eq:valuefun}
& V_N(x,\boldsymbol\pi) \triangleq \mathbf{E}\bigg[ \sum_{i=0}^{N-1} \| \tilde{x}_i \|_Q^2 + \| \pi_i(\tilde{x}_i) \|_R^2 + \| \tilde{x}_N \|^2_Q \bigg];
\end{align}
\textit{$Q \in \mathbb{S}^{n_x}_{+}$ and $R \in \mathbb{S}^{n_u}_{+}$ are positive semidefinite state and input weight matrices, respectively; and $p_w^N \triangleq p_w \times \cdots \times p_w$ denotes the joint PDF of the independent disturbance sequence $\mathbf{w}$ (with $\Delta_w^N \triangleq \Delta_w \times \cdots \times \Delta_w$ being the joint support of $\mathbf{w}$). Let $\boldsymbol\pi^\star$ denote the optimal feedback control policy that minimizes \eqref{e_P1}. At each sampling time, the optimal control inputs $\pi_0^\star$ are applied to the system~\eqref{eq:plantmodel}.}    $\hfill \triangleleft$

Problem 1 is intractable due to: (i) arbitrary form of the state feedback control policy $\boldsymbol\pi$, (ii) unbounded nature of stochastic disturbances, which makes satisfying the hard input constraints impractical, and (iii) nonconvexity and general intractability of the state chance constraints. In addition, guaranteeing feasibility and stability of the SMPC Problem 1 is challenging. In this paper, the latter issues are addressed by introducing various approximations. 

To this end, an affine parametrization is adopted for the feedback control laws $\{\pi_i(\cdot)\}_{i=0}^{N-1}$. Optimizing~\eqref{e_P1} directly over a general class of control laws would result in a dynamic programming problem (e.g., see \cite{lee04}), which is computationally intractable for real-time control. Inspired by \cite{lofberg03,bental04,gou06}, this work uses an affine disturbance feedback parametrization for the control laws to arrive at a suboptimal, but tractable convex optimization problem. Hard input constraints are satisfied by saturating the disturbance terms in the parametrized control laws \cite{chatt11}. The Cantelli-Chebyshev inequality \cite{mar79} is used to replace the state chance constraints with a deterministic surrogate in terms of the mean and variance of the predicted states. Next, the methods employed to obtain a tractable program for Problem 1 are discussed.


\section{Methods for Tractable SMPC Problem}
\label{sec:methods}

\subsection{Compact notation}

Predictions of the dynamics of \eqref{eq:plantmodel} over the prediction horizon $N$ are required to obtain a tractable formulation for Problem 1. Define the stacked vectors $\mathbf{x} \in \mathbb{R}^{n_x(N+1)}$ and $\mathbf{u} \in \mathbb{R}^{n_uN}$ to represent, respectively, the system states (predicted by~\eqref{eq:statepred}) and the control policy over the prediction horizon, i.e.,
\begin{align*}
& \mathbf{x} \triangleq [\tilde{x}_0^\top,\ldots,\tilde{x}_N^\top]^\top \\
& \mathbf{u} \triangleq [\pi_0^\top,\ldots,\pi_{N-1}^\top]^\top
\end{align*}
with $\tilde{x}_0 = x$ being the observed states. The dynamics of \eqref{eq:plantmodel} are compactly described by 
\begin{align} \label{eq:predmodelcompact}
\mathbf{x} = \mathbf{A}\tilde{x}_0 + \mathbf{B}\mathbf{u} + \mathbf{D}\mathbf{G}\mathbf{w},
\end{align}
where
\begin{align*}
\mathbf{A} \triangleq \begin{bmatrix}
I_{n_x} \\
A \\
\vdots \\
A^N
\end{bmatrix};~\mathbf{D} \triangleq \begin{bmatrix}
0 & \cdots & \cdots & 0 \\
I_{n_x} & \ddots &  & \vdots \\
A & I_{n_x} & \ddots & \vdots \\
\vdots & & \ddots & 0 \\
A^{N-1} & \cdots & A & I_{n_x}
\end{bmatrix};
\end{align*}
$\mathbf{B} \triangleq \mathbf{D}(I_N \otimes B)$; and $\mathbf{G} \triangleq I_N \otimes G$. Hence, the value function~\eqref{eq:valuefun} can be written in the compact form
\begin{align} \label{eq:valuefuncompact}
V_N(x,\mathbf{u}) = \mathbf{E}\Big[ \| \mathbf{x} \|^2_\mathbf{Q} + \| \mathbf{u} \|^2_\mathbf{R} \Big],
\end{align}
where $\mathbf{Q} \triangleq \text{diag}(Q,\ldots,Q)$; and $\mathbf{R} \triangleq \text{diag}(R,\ldots,R)$. Define the set $\mathbb{U}^N \triangleq \mathbb{U} \times \cdots \times \mathbb{U}$ so that $\mathbf{u} \in \mathbb{U}^N$ should hold for all disturbance realizations. From Assumption 1(iv), the set $\mathbb{U}^N$ is a polytope defined by
\begin{align*}
\mathbb{U}^N \triangleq \mathbb{U} \times \cdots \times \mathbb{U} = \{ \mathbf{u}\in\mathbb{R}^{n_u N} \mid \mathbf{H}_u \mathbf{u} \leq \mathbf{k}_u \},
\end{align*}
where $\mathbf{H}_u \in \mathbb{R}^{sN \times n_u N}$ and $\mathbf{k}_u \in \mathbb{R}^{sN}$ are given by $\mathbf{H}_u \triangleq \text{diag}(H_u,\ldots,H_u)$ and $\mathbf{k}_u \triangleq [k_u^\top,\ldots,k_u^\top]^\top$, respectively. Likewise, the JCCs in \eqref{e_P1} can be compactly written as
\begin{align} \label{eq:JCCcompact}
\mathbf{Pr}\big[ \mathbf{H}_x^{[i]} \mathbf{x} \leq k_x \big] \geq 1 - \beta,~\forall i \in \mathbb{Z}_{[1,N]}
\end{align}
with $\mathbf{H}_x^{[i]} \in \mathbb{R}^{r \times n_x(N+1)}$ defined by
\begin{align*}
\mathbf{H}_x^{[i]} \triangleq [
0,\cdots,0,~1,~0,\cdots,0
] \otimes H_x.
\end{align*}
where the $1$ is in the $(i+1)^{\mathrm{th}}$ position of the $(N+1)$-length vector. 

\subsection{Tractable feedback control policy} 

Linear-quadratic-Gaussian (LQG) control minimizes the value function \eqref{eq:valuefun} in the absence of input and state constraints \cite{athans71}. The solution to this problem can be obtained analytically and has the form of a linear state feedback control law. Similarly, to obtain a tractable formulation for~\eqref{e_P1}, each element of $\boldsymbol\pi$ can be parametrized as an affine state feedback control law such that
\begin{align} \label{eq:affinestate}
\pi_i \triangleq \sum_{j=0}^i L_{i,j} \tilde{x}_j + g_i,~\forall i \in \mathbb{Z}_{[0,N-1]},
\end{align}
where $L_{i,j} \in \mathbb{R}^{n_u \times n_x}$ and $g_i \in \mathbb{R}^{n_u}$ would be decision variables in the optimization problem. However, the affine state feedback parametrization \eqref{eq:affinestate} would result in nonconvex optimization \cite{gou06}. In this work, the feedback control laws $\pi_i$ in \eqref{eq:input} are parametrized as affine functions of past disturbances 
\begin{align} \label{eq:affinedist}
\pi_i \triangleq \sum_{j=0}^{i-1} M_{i,j} (Gw_j) + v_i,~\forall i \in \mathbb{Z}_{[0,N-1]}
\end{align}
with $M_{i,j} \in \mathbb{R}^{n_u \times n_x}$ and $v_i \in \mathbb{R}^{n_u}$ being the decision variables.\footnote{The parametrization~\eqref{eq:affinedist} has also been used in \cite{oldew08,hokayem09} to formulate convex SMPC problems.} Due to the assumption of perfect state observation, the disturbance realizations at each time instant can be exactly computed by $Gw_i = \tilde{x}_{i+1} - A\tilde{x}_i - B\pi_i$. Hence, the affine state feedback parametrization \eqref{eq:affinestate} and the affine disturbance feedback parametrization \eqref{eq:affinedist} are equivalent; there exists a nonlinear one-to-one mapping between the two parametrizations \cite[Theorem 9]{gou06}. Even though parametrizing the control policy $\boldsymbol\pi$ using \eqref{eq:affinedist} yields suboptimal solutions to Problem 1, it will allow converting \eqref{e_P1} into a tractable convex problem (as shown in Section \ref{sec:SMPC}). 

Using the affine disturbance feedback parametrization \eqref{eq:affinedist}, the control policy $\mathbf{u}$ can be written compactly as
\begin{align} \label{eq:affinedistcompact}
\mathbf{u} = \mathbf{M}\mathbf{G}\mathbf{w} + \mathbf{v},
\end{align}
where the block lower triangular matrix $\mathbf{M} \in \mathbb{R}^{n_u N \times n_xN}$ and stacked vector $\mathbf{v} \in \mathbb{R}^{n_uN}$ are defined by
\begin{align} \label{eq:Mstruct}
\mathbf{M} \triangleq \begin{bmatrix}
0 & \cdots & \cdots & 0 \\
M_{1,0} & 0 &\cdots &0 \\
\vdots &\ddots &\ddots &\vdots \\
M_{N-1,0} &\cdots &M_{N-1,N-2} &0
\end{bmatrix}\; \text{and} \; \mathbf{v}\triangleq\begin{bmatrix}
v_0 \\
v_1 \\
\vdots \\
v_{N-1}
\end{bmatrix},
\end{align}
respectively. Note that $(\mathbf{M},\mathbf{v})$ will comprise the decision variables in the tractable surrogate for Problem 1. Given states $x$ and the control policy \eqref{eq:affinedistcompact}, a pair $(\mathbf{M},\mathbf{v})$ is said to be admissible for Problem 1 if $(\mathbf{M},\mathbf{v})$ results in state and input sequences that satisfy the constraints~\eqref{eq:inputcons} and \eqref{eq:statecc} over the prediction horizon. Hence, the set of admissible decision variables is defined by
\begin{align} \label{eq:PiNx}
\Pi_N(x) \triangleq \left\{
(\mathbf{M},\mathbf{v})\left\vert
\begin{array}{r}
(\mathbf{M},\mathbf{v})\text{ satisfies \eqref{eq:Mstruct}}, \\
\mathbf{x} = \mathbf{A}\tilde{x}_0 + \mathbf{B}\mathbf{u} + \mathbf{D}\mathbf{G}\mathbf{w}, \\
\tilde{x}_0 = x,~\mathbf{u} = \mathbf{M}\mathbf{G}\mathbf{w} + \mathbf{v}, \\
\mathbf{Pr}\big[\mathbf{H}_u\mathbf{u} \leq \mathbf{k}_u\big] = 1, \\
\mathbf{Pr}\big[ \mathbf{H}_{x}^{[i]} \mathbf{x} \leq k_x \big] \geq 1 - \beta, \\
\mathbf{w} \sim p_w^N,~\forall i \in \mathbb{Z}_{[1,N]}
\end{array}\right. \right\}.
\end{align}

\subsection{Saturation of stochastic disturbances for handling hard input constraints}

In general, it is impractical to use a linear state feedback control law for ensuring the hard input constraints \eqref{eq:inputcons} in the presence of unbounded disturbances. Hence, the set $\Pi_N(x)$ can only admit solutions with $\mathbf{M} = 0$ when the joint PDF $p_w^N$ has unbounded support. This can, however, potentially lead to a loss in control performance since the control policy $\mathbf{u} = \mathbf{v}$ will not account for knowledge of future disturbances in the state prediction. To allow for handling hard input constraints in the presence of unbounded disturbances, the input constraints \eqref{eq:inputcons} can be relaxed in terms of expectation-type constraints \cite{primbs09} or chance constraints \cite{farina15}. However, these approaches would not provide a rigorous guarantee for fulfillment of \eqref{eq:inputcons} with respect to all disturbance realizations.\footnote{Note that when the disturbances lie in a compact set (i.e., bounded disturbances), hard input bounds can be guaranteed using worst-case robust control approaches.}  

In this work, the disturbance terms in the control policy \eqref{eq:affinedistcompact} are saturated to enable dealing with unbounded disturbances \cite{hokayem12,chatt11}. The control policy \eqref{eq:affinedistcompact} in the admissible set $\Pi_N(x)$ is replaced with
\begin{align} \label{eq:affinedistcompactsat}
\mathbf{u} = \mathbf{M}\boldsymbol\varphi(\mathbf{G}\mathbf{w}) + \mathbf{v},
\end{align}
where $\boldsymbol\varphi(\mathbf{G}\mathbf{w}) \triangleq [\varphi_0(Gw_0)^\top,\ldots,\varphi_{N-1}(Gw_{N-1})^\top]^\top$; $\varphi_i(Gw_i) \triangleq [\varphi^1_i(G_1 w_i),\ldots,\varphi^{n_x}_i(G_{n_x} w_{i})]^\top$; $G_i$ is the $i^{\mathrm{th}}$ row of the matrix $G$; and $\varphi_i^j : \mathbb{R} \rightarrow \mathbb{R}$ denotes saturation functions with the property $\sup_{a \in \mathbb{R}} | \varphi_i^j(a) | \leq \varphi_{\text{max}}, \, \forall i \in \mathbb{Z}_{[0,N-1]}, \, \forall j \in \mathbb{Z}_{[1,n_x]}$. The latter property implies that $\|\boldsymbol\varphi(\mathbf{G}\mathbf{w}) \|_{\infty} \leq \varphi_{\text{max}}$. The choice of the element-wise saturation functions $\varphi_i^j$ is arbitrary; commonly used saturation functions are piecewise linear or sigmoidal functions \cite{hokayem09}. Since $\|\boldsymbol\varphi(\mathbf{G}\mathbf{w}) \|_{\infty} \leq \varphi_{\text{max}}, \, \forall \mathbf{w} \in \Delta_w^N$, the saturated disturbance constraints are defined in terms of the polytopic set 
\begin{align}
\mathcal{W} \triangleq \{ \mathbf{w} \in \mathbb{R}^{n_x N} \mid \mathbf{H}_w \mathbf{w} \leq \mathbf{k}_w\},
\end{align}
where $\mathbf{H}_w \in \mathbb{R}^{a \times n_x N}$; $\mathbf{k}_w\in\mathbb{R}^{a}$; and $a \in \mathbb{N}$ is the number of saturated disturbance constraints. The input constraint $\mathbf{Pr}\big[\mathbf{H}_u\mathbf{u} \leq \mathbf{k}_u\big] = 1$ associated with the control policy \eqref{eq:affinedistcompactsat} can now be rewritten as
\begin{align} \label{eq:robustinputcons}
\mathbf{H}_u \mathbf{v} + \max_{\boldsymbol\varphi(\mathbf{G}\mathbf{w}) \in \mathcal{W}}(\mathbf{H}_u \mathbf{M}\boldsymbol\varphi(\mathbf{G}\mathbf{w})) \leq \mathbf{k}_u,
\end{align}
where the maximization is row-wise (i.e., maximum of each element in the vector). The following Lemma indicates that \eqref{eq:robustinputcons} can be rewritten in terms of a set of convex inequalities.

\textbf{Lemma 1.} \textit{The input constraint \eqref{eq:robustinputcons} can be represented exactly by linear inequalities $\mathbf{H}_u \mathbf{v} + \mathbf{Z}^\top \mathbf{k}_w \leq \mathbf{k}_u$ and $\mathbf{Z} \geq 0$ (element-wise) for any $\mathbf{Z} \in \mathbb{R}^{a \times sN}$ satisfying $\mathbf{Z}^\top \mathbf{H}_w = \mathbf{H}_u \mathbf{M}$.}

\subsection{Joint chance constraints}

The joint chance constraints~\eqref{eq:statecc} are generally nonconvex even when the affine control policy \eqref{eq:affinedistcompactsat} is used to solve Problem 1. Mutual exclusivity of events $\textbf{Pr}[x \in \mathbb{X}]$ and $\textbf{Pr}[x \not\in \mathbb{X}]$ can be used to obtain a convex relaxation for \eqref{eq:statecc}. First, note that $\textbf{Pr}[x \in \mathbb{X}] + \textbf{Pr}[x \not\in \mathbb{X}] = 1$ such that $\textbf{Pr}[x \not\in \mathbb{X}] \leq \beta$ is equivalent to \eqref{eq:statecc}. Using \eqref{eq:stateconsset}, it can be derived that $\textbf{Pr}[x \not\in \mathbb{X}] = \textbf{Pr}[\cup_{i=1}^{r} (H_{x,i} x \geq k_{x,i})]$, where $H_{x,i}$ and $k_{x,i}$ are the $i^{\mathrm{th}}$ row of $H_x$ and $k_x$, respectively. As the latter term is the union of a set of events, Boole's inequality (aka the union bound) is applied to obtain $\textbf{Pr}[x \not\in \mathbb{X}] \leq \sum_{i=1}^r \textbf{Pr}[H_{x,i} x \geq k_{x,i}]$. Requiring the right-hand side of this inequality to be less than or equal to $\beta$ ensures that \eqref{eq:statecc} holds.

The above result directly extends to \eqref{eq:JCCcompact}. To obtain a convex surrogate for~\eqref{eq:JCCcompact}, each JCC is written in terms of $r$ individual chance constraints (ICCs) of the form
\begin{align} \label{eq:ICCscompact}
\mathbf{Pr}\big[ \mathbf{H}^{[i]}_{x,j} \mathbf{x} \geq k_{x,j} \big] \leq \alpha_j, \quad \forall i \in \mathbb{Z}_{[1,N]}, \, j \in \mathbb{Z}_{[1,r]},  
\end{align}
where $\mathbf{H}^{[i]}_{x,j}$ is the $j^{\mathrm{th}}$ row of the matrix $\mathbf{H}_x^{[i]}$ and $\alpha_j \in [0,1)$ is the maximum probability of violation for the $j^{\mathrm{th}}$ constraint in the polytope~\eqref{eq:stateconsset}. Whenever $\alpha_1 + \alpha_2 + \cdots + \alpha_r = \beta$, the constraints \eqref{eq:JCCcompact} are guaranteed to hold from the analysis above. The Cantelli-Chebyshev inequality is used to convert the ICCs \eqref{eq:ICCscompact} into deterministic, but tighter constraints, in terms of the mean and covariance of the predicted states (see also \cite{farina15}).

\textbf{Lemma 2 (Cantelli-Chebyshev Inequality \cite{mar79}).} \textit{Let $Z$ be a scalar random variable with finite variance. For every $a \in \mathbb{R}_+$, it holds that}
\begin{align} \label{eq:CantelliIC}
\textbf{Pr}\big[Z \geq \mathbf{E}[Z] + a \big] \leq \frac{\mathbf{Var}[Z]}{\mathbf{Var}[Z] + a^2}.
\end{align}

\textbf{Lemma 3.} \textit{Let $x$ be a vector of random variables with mean $\bar{x} \in \mathbb{R}^n$ and finite covariance $\Sigma_x \in \mathbb{S}_{++}^n$, $a \in \mathbb{R}^n$, and $b \in \mathbb{R}$. Then, the ICC}
\begin{align*}
\textbf{Pr}[a^\top x \geq b] \leq \varepsilon
\end{align*}
\textit{is equivalent to the deterministic constraints}
\begin{gather*}
a^\top \bar{x} \leq b - \delta \\
a^\top \Sigma_x a \leq \frac{\varepsilon \delta^2}{1 - \varepsilon}
\end{gather*}
\textit{for any $\delta \in \mathbb{R}_+$ and $\varepsilon \in [0,1)$.}

According to Lemma 3\footnote{A similar result is derived in \cite{farina15}, which is used to constrain the probability of each individual element of the state vector.}, satisfaction of the ICCs~\eqref{eq:ICCscompact} can be guaranteed by imposing
\begin{subequations} \label{eq:ccccccc}
\begin{align} 
&\mathbf{H}_{x,j}^{[i]} \bar{\mathbf{x}} \leq k_{x,j} - \delta_j, &\forall i \in \mathbb{Z}_{[1,N]}, \, j \in \mathbb{Z}_{[1,r]}  \label{eq:expectedvalueterm}  \\ 
&\mathbf{H}^{[i]}_{x,j} \Sigma_{\mathbf{x}} \mathbf{H}^{[i]\top}_{x,j} \leq \frac{\alpha_j \delta_j^2}{1 - \alpha_j},  & \forall i \in \mathbb{Z}_{[1,N]}, \, j \in \mathbb{Z}_{[1,r]}  \label{eq:vartermCC}
\end{align}
\end{subequations}
for any $\delta_j \in \mathbb{R}_+$. A single deterministic inequality can be derived for ICCs~\eqref{eq:ICCscompact} by taking the square root of \eqref{eq:vartermCC} to obtain a lower bound for $\delta_{j}$ and substituting the latter into \eqref{eq:expectedvalueterm}. In this case, the optimizer implicitly searches for a positive $\delta_j$ that satisfies \eqref{eq:expectedvalueterm} and \eqref{eq:vartermCC}. However, these constraints have a nonlinear dependence on $\Sigma_\mathbf{x}$, implying that they are nonconvex. The expression can be convexified by linearization around some state values, as has been done in \cite{farina15}. This approach, however, leads to tightening of the constraints and introduces an additional design parameter (the state around which to linearize) that can be hard to determine \textit{a priori}. In this work, \eqref{eq:expectedvalueterm} and \eqref{eq:vartermCC} are treated individually since they are convex in the decision variables for a fixed $\delta_j$ (as shown in Section \ref{sec:SMPC}). Although fixing $\delta_j$ can introduce additional conservatism, the parameters $\delta_j$ relate directly to the size of the states' covariance such that there is a physical meaning associated with their values.

\textbf{Remark 1}
\textit{As discussed in \cite[Proposition 3]{cin11}, the JCCs in \eqref{e_P1} is convex when \eqref{eq:affinedistcompact} is used with a fixed $\mathbf{M}$ and $p_w$ is log-concave. This notion is a commonly used in SMPC, where control policies involve a fixed prestabilizing feedback in conjunction with online optimization of control actions (see e.g., \cite{cannon09,wan09}). However, such control policies would introduce additional suboptimality in the control law as they do not take feedback into account during prediction in the SMPC problem.} $\hfill \triangleleft$

\section{Proposed Approach for Stochastic Model Predictive Control}
\label{sec:SMPC}

In this section, the methods discussed in Section~\ref{sec:methods} are used to derive a tractable convex formulation for the SMPC Problem 1. Feasibility and stability of the proposed SMPC approach are established.   

\subsection{Convex Formulation for SMPC}

To cast Problem 1 in terms of a deterministic convex optimization problem, the set of admissible decision variables \eqref{eq:PiNx} is adapted based on the saturated affine disturbance feedback control policy \eqref{eq:affinedistcompactsat} and the chance constraint approximations \eqref{eq:ccccccc}. Using the system model~\eqref{eq:predmodelcompact}, the dynamics for the mean and variance of the states $\mathbf{x}$ are described by
\begin{gather}
\label{eq:xmeancompact}
\bar{\mathbf{x}} = \mathbf{A}x_0 + \mathbf{B}\bar{\mathbf{u}} \\[1mm]
\label{eq:Sigmaxcompact}
\Sigma_{\mathbf{x}} = \mathbf{B}\Sigma_{\mathbf{u}}\mathbf{B}^\top + \mathbf{B}\sigma_{\mathbf{u},\mathbf{w}}\mathbf{G}^\top\mathbf{D}^\top \\\notag 
~~~~~~~~~~~~~~ + \mathbf{D}\mathbf{G}\sigma_{\mathbf{u},\mathbf{w}}^\top\mathbf{B}^\top +  \mathbf{D}\mathbf{G}\Sigma_{\mathbf{w}}\mathbf{G}^\top\mathbf{D}^\top.
\end{gather}
Since the state constraints \eqref{eq:ccccccc} are merely a function of $\bar{\mathbf{x}}$ and $\Sigma_{\mathbf{x}}$, the set of admissible decision variables based on the control policy \eqref{eq:affinedistcompactsat} takes the form  
\begin{align} \label{eq:admisset}
\Pi^{d}_N(x) \triangleq \left\{
(\mathbf{M},\mathbf{v})\left\vert
\begin{array}{r}
(\bar{\mathbf{x}},\Sigma_{\mathbf{x}}) \text{ given by \eqref{eq:xmeancompact} and \eqref{eq:Sigmaxcompact}}, \\
(\mathbf{M},\mathbf{v})\text{ satisfies \eqref{eq:Mstruct}},~\tilde{x}_0 = x, \\
\mathbf{u} = \mathbf{M}\boldsymbol\varphi(\mathbf{G}\mathbf{w}) + \mathbf{v}, \\
\mathbf{Pr}\big[\mathbf{H}_u\mathbf{u} \leq \mathbf{k}_u\big] = 1,\\
\mathbf{H}_{x,j}^{[i]} \bar{\mathbf{x}} \leq k_{x,j} - \delta_j, \\
\mathbf{H}^{[i]}_{x,j} \Sigma_{\mathbf{x}} \mathbf{H}^{[i]\top}_{x,j} \leq \alpha_j \delta_j^2(1 - \alpha_j)^{-1}, \\
\forall i \in \mathbb{Z}_{[1,N]},~\forall j \in \mathbb{Z}_{[1,r]},~\mathbf{w} \sim p_w^N 
\end{array}\right. \right\}
\end{align}
for any $\delta_j \in \mathbb{R}_+$ and $\alpha_j \in [0,1)$ such that $\alpha_1 + \alpha_2 + \cdots + \alpha_r = \beta$. The SMPC Problem 1 can now be rewritten as the convex stochastic optimal control problem\footnote{An explicit expression for $V_N(x,\mathbf{M},\mathbf{v})$ is derived in the proof of Theorem 1.}
\begin{align} \label{eq:P2}
(\mathbf{M}^\star(x),\mathbf{v}^\star(x)) \triangleq \argmin_{(\mathbf{M},\mathbf{v}) \in \Pi^{d}_N(x)} V_N(x,\mathbf{M},\mathbf{v}),
\end{align}
where the optimal control policy \eqref{eq:affinedistcompactsat} with $(\mathbf{M},\mathbf{v}) = (\mathbf{M}^\star(x),\mathbf{v}^\star(x))$, if one exists, will generate states that satisfy the hard input constraints \eqref{eq:inputcons} and the joint state chance constraints \eqref{eq:statecc}. The set of feasible initial conditions, for which a feasible controller of the form \eqref{eq:affinedistcompactsat} exists in \eqref{eq:P2} (i.e., the domain of attraction of \eqref{eq:P2}), is defined by
\begin{align} \label{eq:domainofattract}
\mathcal{X}^d_N \triangleq \left\{ x \in \mathbb{R}^{n_x} \mid \Pi^{d}_N(x) \neq \emptyset \right\}.
\end{align}
Note that the domain of attraction $\mathcal{X}^d_N$ is a function of the parameters $\delta_j$ and $\alpha_j$ (see \eqref{eq:ccccccc}), which are prespecified in \eqref{eq:P2}.

\subsection{Feasibility of the Convex SMPC Formulation}

The inclusion of the state constraints \eqref{eq:ccccccc} within \eqref{eq:admisset} renders the domain of attraction $\mathcal{X}^d_N$ to be a subset of $\mathbb{R}^{n_x}$ (i.e., $\mathcal{X}^d_N \subset \mathbb{R}^{n_x}, \, \forall N \in \mathbb{N}$) as, for certain observed states $x$, there will not exist any controller for which future state predictions satisfy the constraints. In addition, when the stochastic disturbances $w$ are unbounded, the states of the true system \eqref{eq:plantmodel} will violate any given compact set infinitely often over an infinite time horizon for any bounded control action $u \in \mathbb{U}$. This implies that it cannot be guaranteed that $x^+ \in \mathcal{X}^d_N$ even when $x \in \mathcal{X}^d_N$. Hence, it is impossible to ensure feasibility of the SMPC problem \eqref{eq:P2} at all times due to the unbounded nature of the stochastic disturbances. 

To guarantee feasibility of the stochastic optimal control problem~\eqref{eq:P2}, this work considers \textsl{softening} the state constraints~\eqref{eq:ccccccc} to ensure the program is always feasible and results in minimal constraint violation \cite{ker00}.\footnote{Note that the underlying notion of the chance constraints~\eqref{eq:statecc} inherently involves softening of the hard state constraints $x \in \mathbb{X}$ for any possible disturbances.} A systematic way to soften the constraints \eqref{eq:ccccccc} is to introduce slack variables $\boldsymbol\epsilon$ in the optimal control problem \eqref{eq:P2}, where the magnitude of the slack variables $\boldsymbol\epsilon$ corresponds to the associated amount of constraint violation. A penalty on $\boldsymbol\epsilon$ is then included in the value function to minimize the constraint violation, along with minimizing the original cost function. Note that it is desired for the soft-constrained optimal control problem to yield the same solution as the original hard-constrained problem when the latter is feasible. This can be achieved using the exact penalty function method (e.g., see \cite{ker00,hov01,zei10,hov11} and the references therein for details).

\textbf{Theorem 1.} \textit{Let Assumption 1 hold and assume that $\mathbf{E}[\boldsymbol\varphi(\mathbf{G}\mathbf{w})] = 0$. Then, the stochastic optimal control problem \eqref{eq:P2} with softened state constraints is a convex Second Order Cone Program (SOCP) with respect to decision variables $(\mathbf{M},\mathbf{v},\mathbf{Z},\boldsymbol\epsilon)$, and is defined by}
\begin{subequations} \label{e_P3}
\begin{align}
\min_{(\mathbf{M},\mathbf{v},\mathbf{Z},\boldsymbol\epsilon)} ~ \mathbf{b}^\top\mathbf{v} + \|\mathbf{v}\|_{\mathbf{S}_1}^2 + \nu^\top \mathbf{m} + \|\mathbf{m}\|_\Lambda^2 + \rho\mathbf{1}^\top\boldsymbol\epsilon
\end{align}
\end{subequations}
\begin{subequations}
\begin{alignat*}{2}
\text{s.t.: } ~ & \tilde{x}_0  &&= x \tag{\ref{e_P3}b}\\
& (\mathbf{M},\mathbf{v}) &&\text{ satisfies \eqref{eq:Mstruct}} \tag{\ref{e_P3}c} \\
& \mathbf{m} &&= \text{vec}(\mathbf{M}) \tag{\ref{e_P3}d}\\
& \mathbf{H}_u\mathbf{v} + \mathbf{Z}^\top\mathbf{k}_w && \leq \mathbf{k}_u \tag{\ref{e_P3}e}\\
& \mathbf{Z}^\top \mathbf{H}_w &&= \mathbf{H}_u\mathbf{M} \tag{\ref{e_P3}f}\\
& \mathbf{Z}  &&\geq 0 \tag{\ref{e_P3}g}\\
& \mathbf{H}^{[i]}_{x,j} \mathbf{A}\tilde{x}_0 + \mathbf{H}^{[i]}_{x,j}\mathbf{B}\mathbf{v} &&\leq k_{x,j} - \delta_j + \epsilon^m_{i,j}, \label{eq:softconsexp}\tag{\ref{e_P3}h}\\
& \mathbf{y}^\top_{i,j} &&= \mathbf{H}^{[i]}_{x,j}\mathbf{B}\mathbf{M},\tag{\ref{e_P3}i} \\
& \|\mathbf{y}_{i,j}\|_{\Omega_1}^2 + \mathbf{q}_{i,j}^\top \mathbf{y}_{i,j} + c_{i,j} &&\leq \epsilon^v_{i,j}, \tag{\ref{e_P3}j} \label{eq:softconsvar}\\
& \epsilon^m_{i,j} &&\geq 0, \tag{\ref{e_P3}k}\\
& \epsilon^v_{i,j} &&\geq 0, \tag{\ref{e_P3}l}\\
& \forall i \in \mathbb{Z}_{[1,N]},~\forall j \in \mathbb{Z}_{[1,r]},
\end{alignat*}
\end{subequations}
\textit{where}
\begin{align*}
\begin{array}{ll}
\mathbf{b}^\top & \triangleq 2\tilde{x}_0^\top \mathbf{A}^\top\mathbf{Q}\mathbf{B};  \\
\mathbf{S}_1 & \triangleq \mathbf{B}^\top\mathbf{Q}\mathbf{B}+\mathbf{R}; \\
\mathbf{S}_2 & \triangleq 2\mathbf{G}^\top\mathbf{D}^\top\mathbf{Q}\mathbf{B}; \\
\Omega_1 & \triangleq \mathbf{E}[\boldsymbol\varphi(\mathbf{G}\mathbf{w})\boldsymbol\varphi(\mathbf{G}\mathbf{w})^\top]; \\
\Omega_2 & \triangleq \mathbf{E}[\boldsymbol\varphi(\mathbf{G}\mathbf{w})\mathbf{w}^\top]; \\
\Lambda & \triangleq \Omega_1 \otimes S_1^\top; \\
\nu & \triangleq \text{vec}(\mathbf{S}_2^\top\Omega_2);\\
\mathbf{q}_{i,j}^\top & \triangleq 2\mathbf{H}^{[i]}_{x,j}\mathbf{D}\mathbf{G}\Omega_2; \\
c_{i,j} & \triangleq \mathbf{H}^{[i]}_{x,j}\mathbf{D}\mathbf{G}\Sigma_\mathbf{w}\mathbf{G}^\top\mathbf{D}\mathbf{H}^{[i]\top}_{x,j} - \alpha_j\delta_j^2(1-\alpha_j)^{-1};
\end{array}
\end{align*}
\textit{$\boldsymbol\epsilon \in \mathbb{R}^{2rN}$ denotes the vector of all slack variables $\epsilon^m_{i,j}$ and $\epsilon^v_{i,j}$; and $\rho \in \mathbb{R}_+$ is the weight associated with the size of constraint violation.}

Denote the SOCP~\eqref{e_P3} by $\mathbb{P}_N^s(x)$. Define the set of feasible initial conditions for \eqref{e_P3} by
\begin{align}
\mathcal{X}^s_N \triangleq \left\{ x \in \mathbb{R}^{n_x} \mid \exists (\mathbf{M},\mathbf{v},\mathbf{Z},\boldsymbol\epsilon) \text{ feasible to } \mathbb{P}_N^s(x) \right\}.
\end{align}
Since the constraints~\eqref{eq:ccccccc} have been softened in~\eqref{e_P3}, the domain of attraction for $\mathbb{P}_N^s(x)$ is $\mathcal{X}^s_N = \mathbb{R}^{n_x}$. This is because for any $x \in \mathbb{R}^{n_x}$ there will always exist large enough slack variables $\tilde{\boldsymbol\epsilon} \geq 0$ such that $(\mathbf{M},\mathbf{v},\mathbf{Z},\boldsymbol\epsilon) = (0,0,0,\tilde{\boldsymbol\epsilon})$ is a feasible solution to $\mathbb{P}_N^s(x)$ (note that by assumption $0 \in \mathbb{U}^N$).

Let $(\mathbf{M}^{s\star}(x),\mathbf{v}^{s\star}(x))$ be the optimizer of $\mathbb{P}_N^s(x)$. Define the optimal control law $\kappa^s_N : \mathbb{R}^{n_x} \rightarrow \mathbb{U}$ in terms of the saturated affine disturbance feedback parametrization~\eqref{eq:affinedistcompactsat} with $(\mathbf{M},\mathbf{v})=(\mathbf{M}^{s\star}(x),\mathbf{v}^{s\star}(x))$. Receding-horizon control of~\eqref{eq:plantmodel} entails applying the first elements of this optimal control policy, i.e.,        
\begin{align}
\kappa^s_N(x) \triangleq v_0^{s\star}(x),
\end{align}
to \eqref{eq:plantmodel} such that the closed-loop response of the system is
\begin{align} \label{eq:CLsystem}
x^+ = A x + B\kappa^s_N(x) + Gw.
\end{align}
The choice of the weight $\rho$ is critical in solving $\mathbb{P}_N^s(x)$. When the weight $\rho$ is sufficiently large, constraint violations will only occur if there is no feasible solution with hard constraints \cite{de94}. This is known as the exact penalty function method. A condition on the lower bound for $\rho$ is given by \cite{ker00}
\begin{align} \label{eq:rho}
\rho > & \max_{x,\boldsymbol\lambda} \| \boldsymbol\lambda \|_\infty \\\notag
& \begin{array}{ll}
\text{subject to:} & x \in \mathcal{X}^d_N,~\text{KKT conditions for } \mathbb{P}_N^s(x),
\end{array}
\end{align}
where $\boldsymbol\lambda$ denotes the Lagrange multipliers for $\mathbb{P}_N^s(x)$. However, determining $\rho$ using~\eqref{eq:rho} is not straightforward as an explicit characterization of $\mathcal{X}^d_N$ is not always available and the KKT conditions are nonconvex. Alternatively, $\mathbf{M}$ can be set equal to zero in \eqref{e_P3} such that $\mathbb{P}_N^s(x)$ reduces to a quadratic program that will be feasible for any initial conditions as long as $c_{i,j} \leq 0$. A conservative lower bound on $\rho$ can then be obtained by solving a finite number of linear programs (e.g., see \cite{ker00}). 

Another practical method for \textit{a priori} selecting $\rho$ is to choose $\rho$ as large as possible such that numerical issues do not occur when solving \eqref{e_P3}. The idea is that at some point $\rho_{\text{max}} < \infty$, the weight will be so large that changes in the original value function are not numerically discernible in the optimizer. Note that locating an approximate value for $\rho_{\text{max}}$ may be easier than solving \eqref{eq:rho}. 


\subsection{Stability of the Convex SMPC Formulation}

To analyze the stability of the closed-loop system \eqref{eq:CLsystem}, the time of occurrence of the states $x$ and disturbances $w$ should be denoted. This is done with a subscript $t$ such that the closed-loop dynamics \eqref{eq:CLsystem} take the form
\begin{align} \label{eq:CLrecursion}
x_{t+1} = A x_t + B\kappa^s_N(x_t) + Gw_t,~\forall t \in \mathbb{N}_0
\end{align}
for some given initial conditions $x_0$. The closed-loop states $\{ x_t \}_{t \in \mathbb{N}_0}$ generated by \eqref{eq:CLrecursion} represent a discrete-time Markov process as the probability distribution of future states $\{ x_s \}_{s\in\mathbb{Z}_{[t+1,\infty]}}$ is conditionally independent of past states $\{ x_s \}_{s \in \mathbb{Z}_{[0,t-1]}}$ given the current state $x_t$. The stability of Markov processes commonly deals with boundedness of sequences of norm-like functions, e.g., $\{ \mathbf{E}_{x_0}[\| x_t \|_p] \}_{t \in \mathbb{N}_0}$ \cite{mey09}. The stability theory for Markov processes has also been used in the context of SMPC \cite{chatt15,paulson15}. Stochastic stability conditions typically involve a \textit{negative drift condition}. In this work, a \textit{geometric drift condition} is used to establish the stability of SOCP~\eqref{e_P3}, as summarized in the sequel.

\textbf{Lemma 4.} \textit{Let $\{ x_t \}_{t \in \mathbb{N}_0}$ be a Markov process.} \textit{Suppose there exists a function} $V : \mathbb{R}^{n_x} \rightarrow \mathbb{R}_+$\textit{, a compact set }$\mathbb{D} \subset \mathbb{R}^{n_x}$\textit{, and constants }$b \in \mathbb{R}_+$\textit{ and }$\lambda \in [0,1)$\textit{ such that }$\mathbf{E}_{x_0}[V(x_1)] \leq \lambda V(x_0)$\textit{ for all }$x_0 \not\in \mathbb{D}$\textit{ and }$\sup_{x\in\mathbb{D}}\mathbf{E}_{x_0}[V(x_1)] = b$. \textit{Then, the sequence $\{ \mathbf{E}_{x_0}[V(x_t)] \}_{t \in \mathbb{N}_0}$ is bounded for all $x_0 \in \mathbb{R}^{n_x}$.}

Proving stability of the closed-loop system \eqref{eq:CLrecursion} (in a stochastic sense) requires verifying the premises of Lemma 4. In fact, Lemma 4 implies that a \textit{geometric drift condition} is satisfied for all states outside of the compact set $\mathbb{D}$, i.e., $\mathbf{E}_{x_0}[V(x_1)] - V(x_0) \leq -(1-\lambda) V(x_0)$ for all $x_0 \not\in \mathbb{D}$.

\textbf{Theorem 2.} \textit{Under the assumptions of Theorem 1, the closed-loop system \eqref{eq:CLrecursion} is stochastically stable. More precisely, $\sup_{t \in \mathbb{N}_0}\mathbf{E}_{x_0}[\| x_t \|^2] < \infty$ for all $x_0 \in \mathbb{R}^{n_x}$.}

As a result of Theorem 2, the control law $\kappa^s(x)$, which is always feasible, leads to a stochastically stable closed-loop system \eqref{eq:CLrecursion} that has bounded variance for all time. 

\section{Case Study: Stochastic Optimal Control of a Continuous ABE Fermentation Process}
\label{sec:case_study}

The SMPC approach is demonstrated for continuous acetone-butanol-ethanol (ABE) fermentation, which is used for bioconversion of lignocellulosic-derived sugars to high-value added \textsl{drop-in} biofuels \cite{jan12,qur08}. The key dynamical characteristics of continuous ABE fermentation include highly nonlinear dynamics, relatively large number of states corresponding to concentration of different species in the metabolic pathway, and inherent biological stochasticity of the pathway. In this work, the dynamic model presented in \cite{haus11} for continuous ABE fermentation using \textit{Clostridium acetobutylicum} is linearized around a desired steady state that corresponds to the solventogenesis stage of cells (see \cite{garcia11} for metabolic pathways in \textit{clostridia}). The linearized system dynamics are described by \eqref{eq:plantmodel}, which consists of $12$ states and $2$ inputs that are the dilution rate ($D$) of the bioreactor and the inlet glucose concentration ($G_0$). All system states are perturbed by zero-mean white noise with variance $10^{-4}$ mM (i.e., $w\sim\mathcal{N}(0,10^{-4}I_{12})$), which has been inferred from experimental data \cite{lee08,zhe13}. The description of system model is given in Appendix~B.
 
The SMPC Problem~\eqref{e_P1} is formulated in terms of setpoint tracking for the ABE products while satisfying hard input constraints and individual chance constraints on the acidic species (i.e., acetate and butyrate). Table~\ref{T1} lists the settings of the corresponding SOCP.\footnote{The concentrations of acetate and butyrate are constrained within $\pm6\%$ of their steady-state values. The hard input constraints are defined based on industrially relevant values \cite{haus11,ban12}.} Receding-horizon control of the continuous ABE fermentation process involves solving the SOCP~\eqref{e_P3} at every sampling instant that the true system states are observed. The tractable SOCP is solved using the CVX package with the Mosek solver \cite{cvx,gb08}.     

\begin{table}[t!]
	\caption{Settings of the SOCP~\eqref{e_P3} for the continuous ABE fermentation process.}
	\begin{center} 
	\resizebox{8.5cm}{!}{%
		\begin{tabular}{ll}
			\hline \hline  \\
			Sampling time for receding-horizon control & 5 min\\
			N & 10  \\
			$Q$ & \text{diag}(0,10,10,0,0,0,0,0,10,0,0,0) \\
			$R$ & \text{diag}(1,1) \\
			$\alpha_j$ in \eqref{eq:ccccccc} & 0.1 \\ [1.0ex]  \hline
			Setpoints ($t\geq 10$ \text{hr}) & $\begin{array}{ll}
			45.7 &\text{\, mM (Acetone)}\\
			54.4 &\text{\, mM (Butanol)}\\
			7.72 &\text{\, mM (Ethanol)}\\
			\end{array}$ \\ \hline
			State constraints & $\begin{array}{l}
			13.94 \text{ mM} \leq x_{Acetate} \leq 15.72 \text{ mM} \\
			10.86 \text{ mM} \leq x_{Butyrate} \leq 12.25 \text{ mM}\\
			\end{array}$ \\ \hline
			Hard input constraints & $\begin{array}{l}
      0.05 \text{ hr}^{-1} \leq D \leq 0.10 \text { hr}^{-1}\\
			56 \text{ mM} \leq G_0 \leq 389 \text{ mM} \\
			\end{array}$ \\ \hline \hline
		\end{tabular} }
	\end{center}
	\label{T1}
\end{table}

\begin{figure}[t!] 
\centering
\subfigure[Acetone concentration profiles]{
\includegraphics[width=225pt]{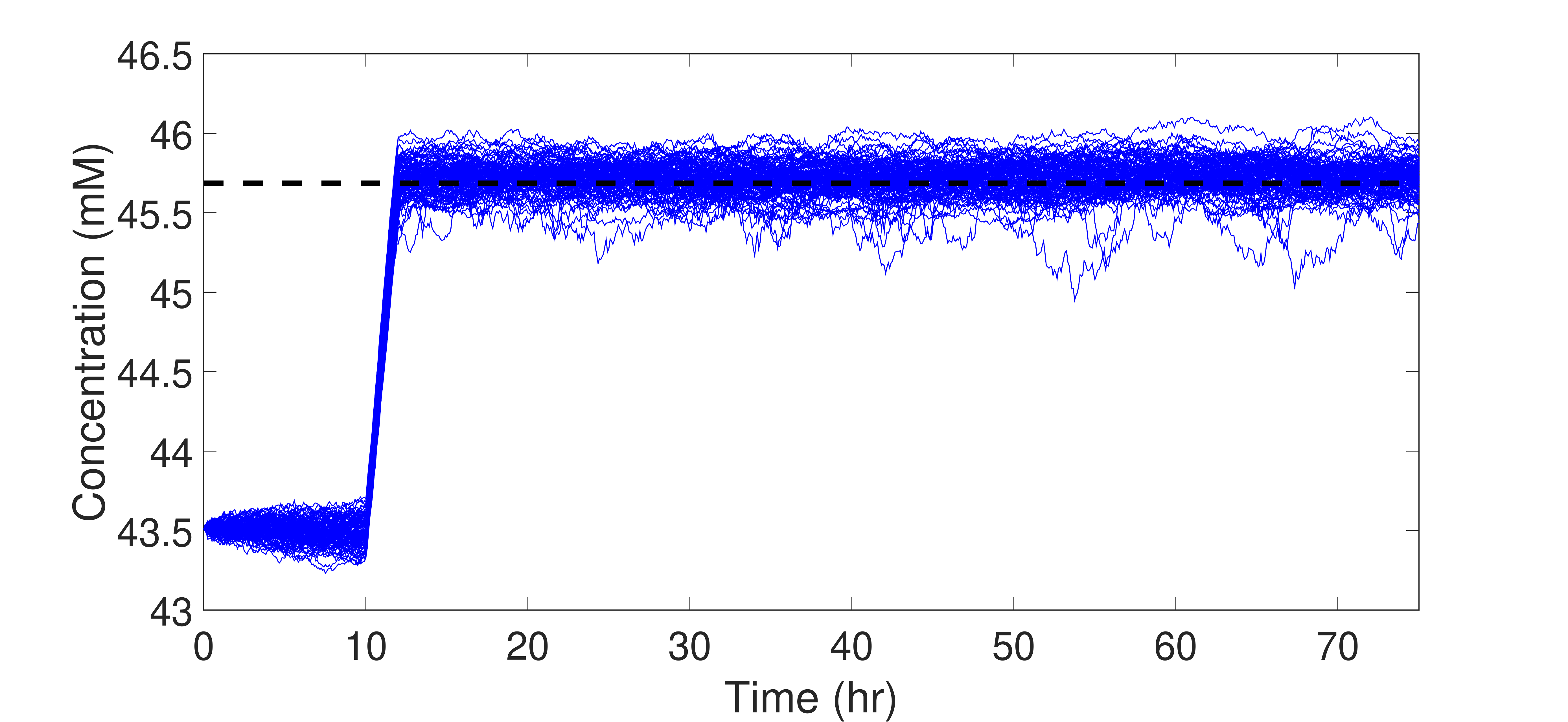}
}
\subfigure[Butanol concentration profiles]{
\includegraphics[width=225pt]{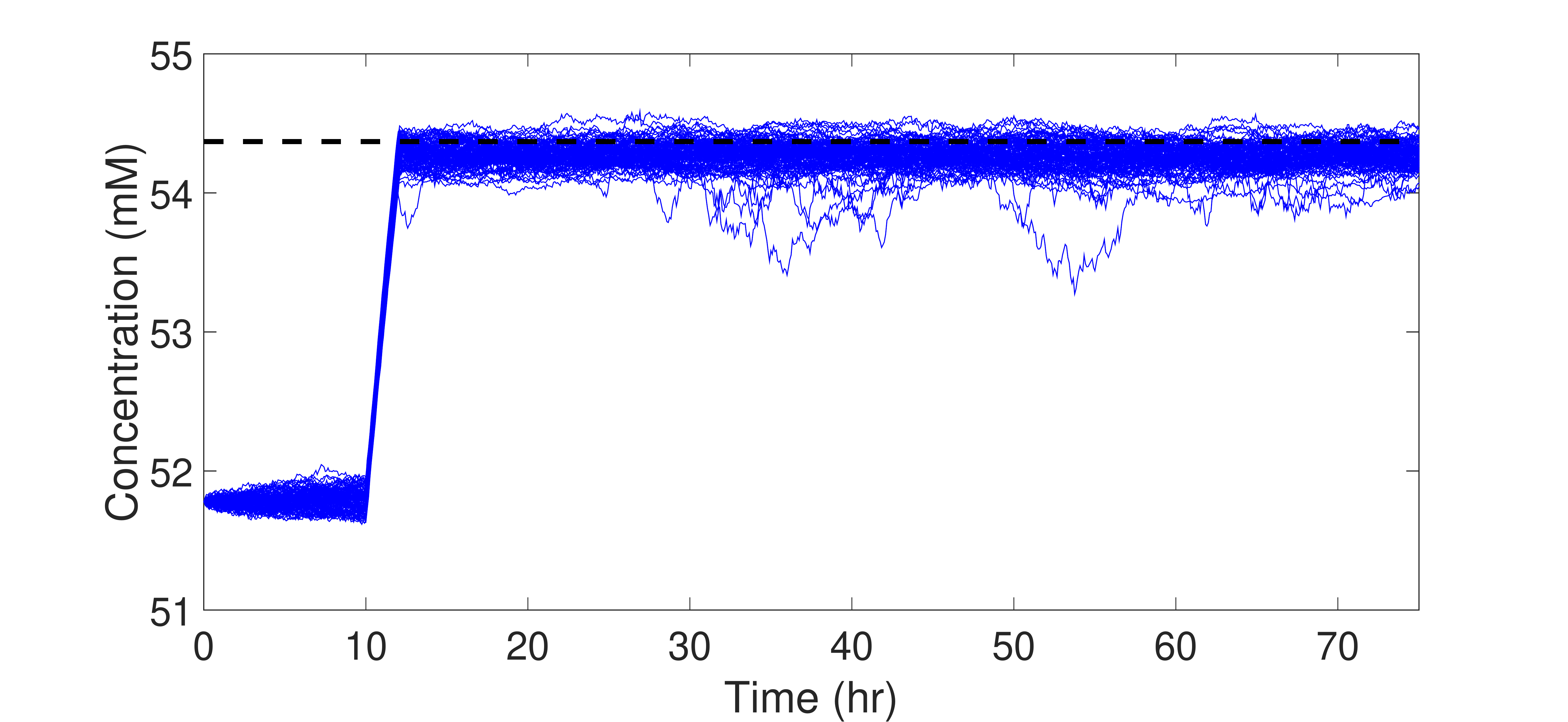}
}
\subfigure[Ethanol concentration profiles]{
\includegraphics[width=225pt]{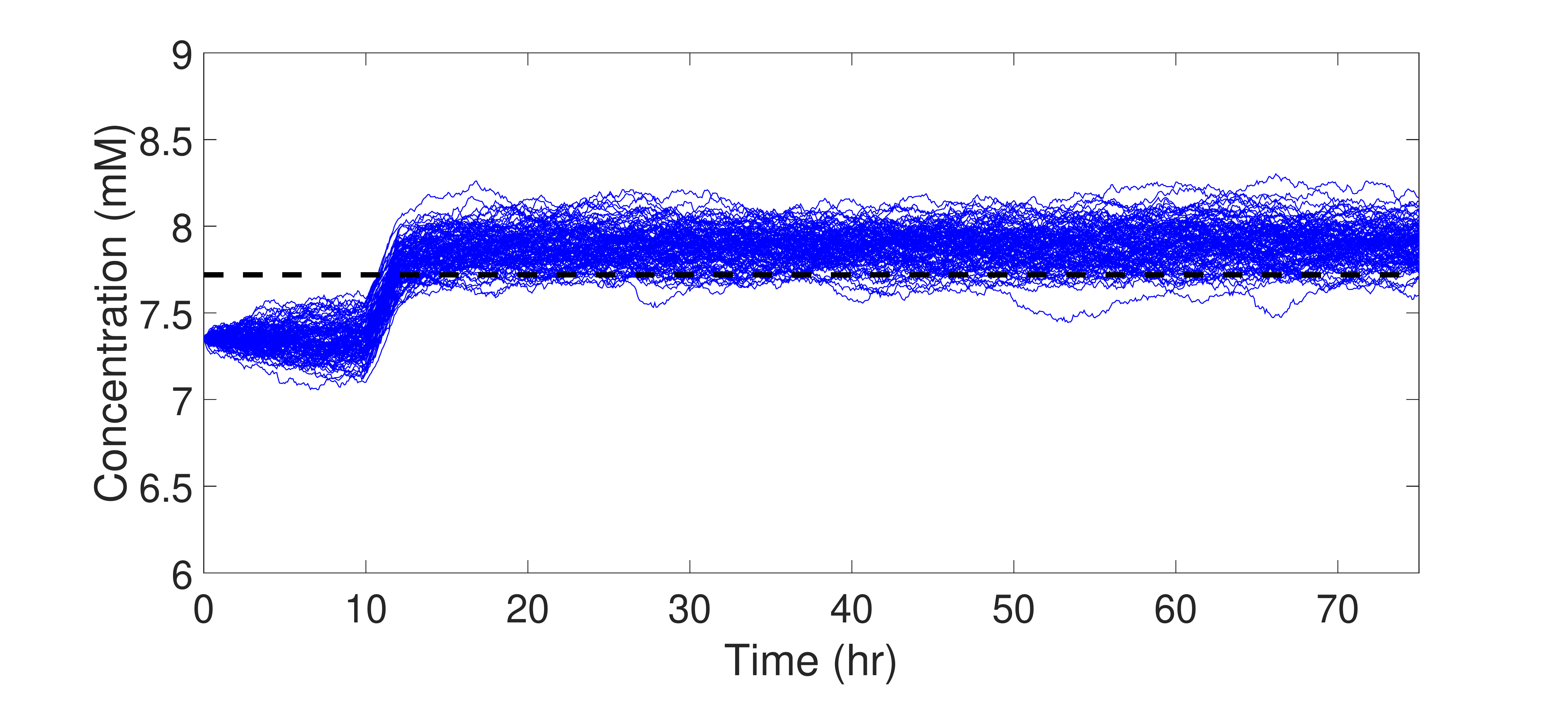}
}
\caption{Concentration profiles of the acetone-butanol-ethanol products in the continuous ABE fermentation process under closed-loop control with the SMPC approach. The concentration profiles are depicted for $100$ closed-loop simulation runs. A $5\%$ change in the setpoints is applied at time $10$ hr.}
\label{fig1}
\end{figure} 

To evaluate the performance of the SMPC approach, $100$ closed-loop simulations have been performed in which the optimal inputs computed by solving the SOCP (i.e., $\kappa_N^s(x)$) at each sampling time instant are applied to the continuous ABE fermentation process. At time $10$ hours, the setpoints for acetone, butanol, and ethanol undergo a $5\%$ increase with respect to their nominal steady-state values. Fig.~\ref{fig1} shows the concentration profiles of acetone, butanol, and ethanol products for the $100$ simulation runs. The SMPC approach allows for effective tracking of the setpoint change, while minimizing the variance of concentration profiles around their respective setpoints. Note that there is a slight offset in setpoint tracking for the product concentration profiles, in particular for the ethanol concentration profiles. The offset results from the tradeoff that the SMPC approach seeks between fulfilling the control objectives and satisfying the state constraints in the presence of process stochasticity (i.e., state chance constraints).
  
\begin{figure}[t!] 
\centering
\subfigure[Probability distributions of acetate concentration]{
\includegraphics[width=225pt]{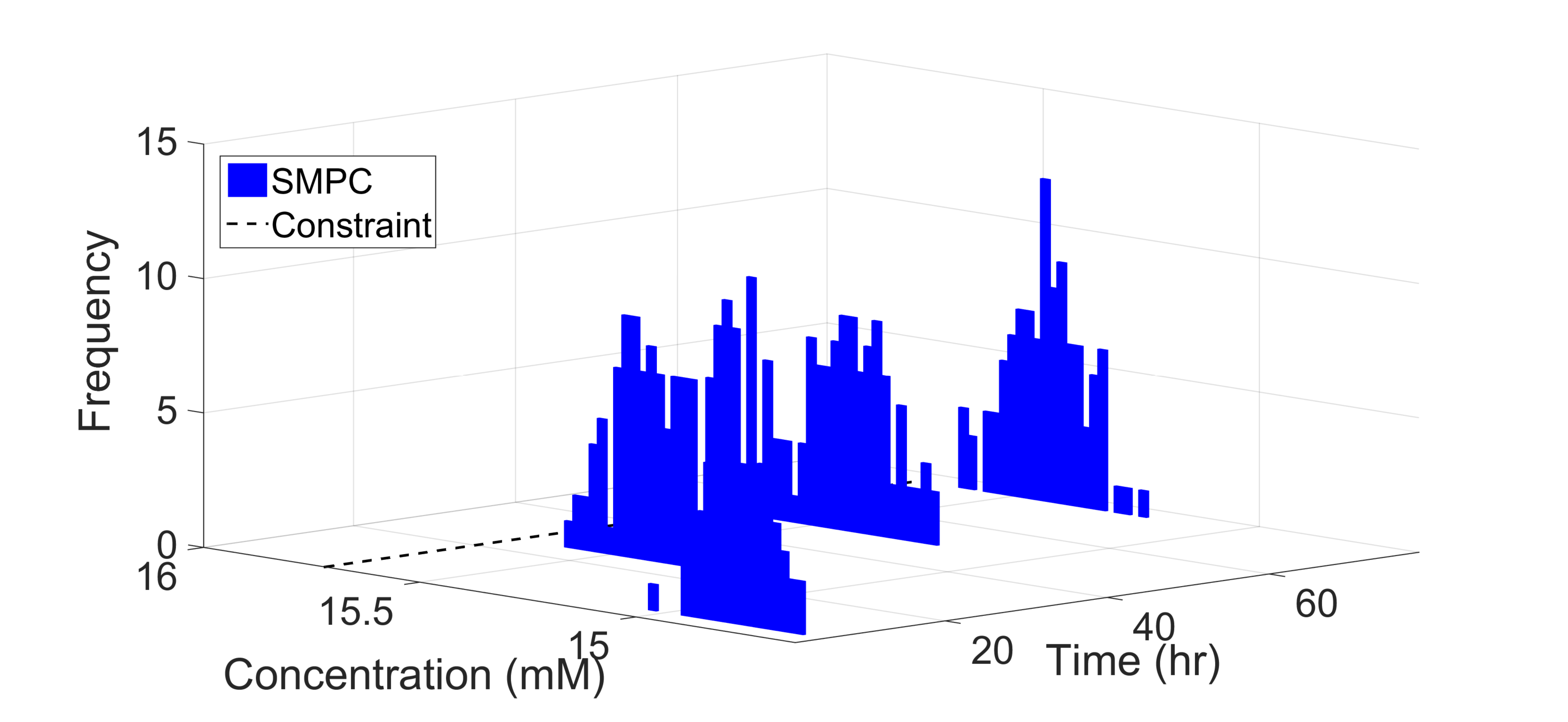}
}
\subfigure[Probability distributions of butyrate concentration]{
\includegraphics[width=225pt]{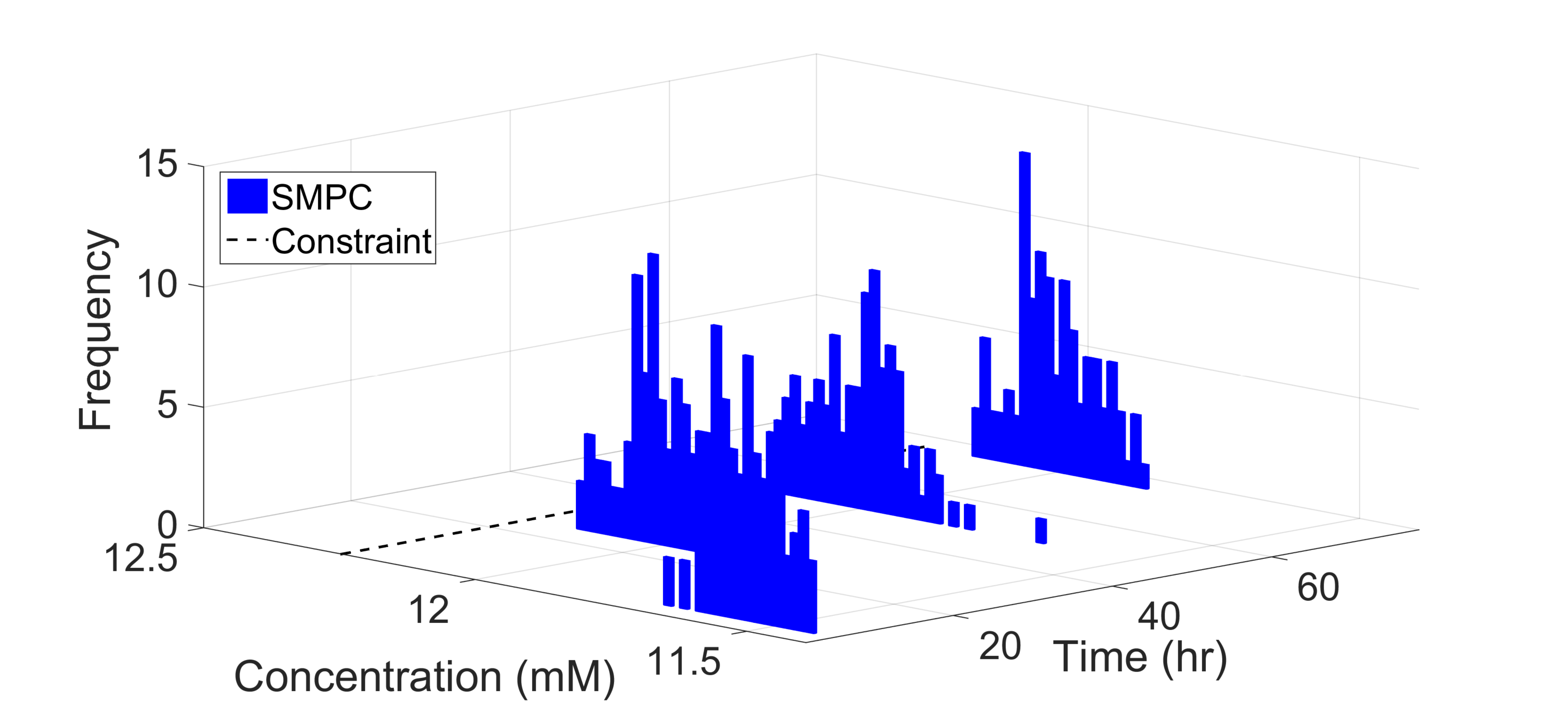}
}
\caption{Evolution of probability distributions of concentration of the acidic species in the continuous ABE fermentation process under closed-loop control with the SMPC approach. The probability distributions are obtained based on $100$ closed-loop simulation runs, and are depicted at four different times. The state constraints are always satisfied using the SMPC approach, whereas the standard MPC controller leads to violation of state constraints in approximately $20\%$ of the closed-loop simulations (not shown here).}
\label{fig2}
\end{figure}

Fig.~\ref{fig2} shows the evolution of probability distributions of the acidic species under SMPC. The state constraints are never violated using SMPC despite allowing for $10\%$ constraint violation. This is attributed to conservative approximation of the chance constraints using the Cantelli-Chebyshev inequality. Seemingly, the SMPC approach applies a hard bound for acetate and butyrate that is more conservative than the bounds considered in the state chance constraints. The ability of SMPC to handle the state constraints in the presence of stochastic uncertainties has been compared to that of a standard MPC controller. The closed-loop simulation results reveal that the standard MPC controller violates the state constraints in approximately $20\%$ of simulations. Maintaining these state constraints is critical to optimal operation of the continuous ABE fermentation, as high concentration of the acidic species can cause the cells to switch from the desired solventogenesis stage to the undesired acidogenesis stage. The simulation case study indicates the capability of the SMPC approach in guaranteeing the fulfillment of state constraints in the presence of process stochasticity at the expense of a slight drop in the control performance.        

\section{Conclusions}
\label{sec:conclusions}

A tractable convex optimization program is derived to solve receding-horizon stochastic model predictive control of linear systems with (possibly) unbounded stochastic disturbances. Due to the unbounded nature of stochastic disturbances, feasibility of the optimization problem cannot always be guaranteed when chance constraints are included in the formulation. To avoid this issue, the state chance constraints are softened to ensure the feasibility. The exact penalty function method is utilized to guarantee the soft-constrained optimization yields the same solution as the hard-constrained problem when the latter is feasible. Building upon the results of \cite{hokayem09}, it is also shown that the receding-horizon control policy results in a stochastically stable system by proving the states satisfy a geometric drift condition outside of a compact set (i.e., the variance of the states is bounded for all time). Future work includes incorporating output feedback with measurement noise, as well as reducing the conservatism introduced by the Cantelli-Chebyshev inequality. 



\section*{Appendix A}
\label{sec:app_A}

This appendix contains the proofs for all Lemmas and Theorems presented in the paper.

\textit{Proof of Lemma 1.} Let the $i^\mathrm{{th}}$ row of the maximization in \eqref{eq:robustinputcons} be the primal linear program (LP). The corresponding dual LP is
\begin{align*}
& \min_{\mathbf{z}_i \geq 0} ~ \mathbf{k}_w^\top \mathbf{z}_i,~\text{s.t.:}~\mathbf{H}_w^\top\mathbf{z}_i = (\mathbf{H}_u \mathbf{M})_i^\top,
\end{align*}
where $(\mathbf{H}_u \mathbf{M})_i$ denotes the $i^{\mathrm{th}}$ row of $\mathbf{H}_u\mathbf{M}$; and $\mathbf{z}_i \in \mathbb{R}^a$ denotes the dual variables. By the weak duality of LPs (e.g., see \cite{boy04}), it is known that the program  $\max_{\boldsymbol\varphi(\mathbf{G}\mathbf{w})\in\mathcal{W}}(\mathbf{H}_u \mathbf{M})_i\boldsymbol\varphi(\mathbf{G}\mathbf{w}) \leq \mathbf{k}_w^\top\mathbf{z}_i$ for any $\mathbf{z}_i \geq 0$ satisfying $\mathbf{H}_w^\top\mathbf{z}_i = (\mathbf{H}_u \mathbf{M})_i^\top$. Stacking the dual variables into the matrix $\mathbf{Z} \triangleq [\mathbf{z}_1,\ldots,\mathbf{z}_{sN}]$ and taking the transpose of this relationship for all $i \in \mathbb{Z}_{[1,sN]}$ results in the inequality $\max_{\boldsymbol\varphi(\mathbf{G}\mathbf{w}) \in \mathcal{W}}(\mathbf{H}_u \mathbf{M}\boldsymbol\varphi(\mathbf{G}\mathbf{w})) \leq \mathbf{Z}^\top \mathbf{k}_w$ for any $\mathbf{Z} \geq 0$ satisfying $\mathbf{Z}^\top \mathbf{H}_w = \mathbf{H}_u\mathbf{M}$. Hence, the assertion of the lemma follows from this statement. $\hfill \square$

\textit{Proof of Lemma 2.} Let $p_Z(z)$ denote the PDF of random variable $Z$ with mean $\bar{Z}$. Define a corresponding zero-mean random variable $Y \triangleq Z - \bar{Z}$ with PDF $p_Y(y)$. Consider $\mathbf{Pr}[Z \geq \bar{Z} + a] = \mathbf{Pr}[Y \geq a] = \int_{a}^\infty p_Y(Y)dz = \mathbf{E}\big[ \mathbf{1}_{[a,\infty)}(Y) \big]$. Define the function $h(y) \triangleq (a+b)^{-2}(y+b)^2$ for any $a \in \mathbb{R}_+$ and $b \in \mathbb{R}_{+}$. It is evident that $\mathbf{1}_{[a,\infty)}(y) \leq h(y), \, \forall y \in \mathbb{R}$ such that $\mathbf{E}\big[ \mathbf{1}_{[a,\infty)}(Y) \big] \leq \mathbf{E}[h(Y)] = (a + b)^{-2}\mathbf{E}[(Y+b)^2] = (a + b)^{-2}(\mathbf{Var}[Y]+b^2)$. The smallest upper bound for $\mathbf{Pr}[Z \geq \bar{Z} + \alpha]$ is obtained by minimizing the right-hand side of the latter inequality with respect to $b$. The solution to this minimization is $b^\star = \mathbf{Var}[Y]/a$. The assertion of Lemma follows by noticing $\mathbf{Var}[Y] = \mathbf{Var}[Z]$ and substituting $b^\star$ for $b$ into this upper bound of $\mathbf{Pr}[Z \geq \bar{Z} + \alpha]$. $\hfill \square$

\textit{Proof of Lemma 3.} When the inequality $a^\top \bar{x} \leq b - \delta$ holds, it is true that
\begin{align*}
\textbf{Pr}[a^\top x \geq b] & \leq \textbf{Pr}[a^\top x \geq a^\top \bar{x} + \delta] \\
& \leq \frac{a^\top \Sigma_x a}{a^\top \Sigma_x a + \delta^2},
\end{align*}
where the second line follows from Lemma 2 for any $\delta\in\mathbb{R}_+$. As the latter inequality is an upper bound for $\textbf{Pr}[a^\top x \geq b] \leq \varepsilon$, the ICC will be satisfied when
\begin{align*}
\frac{a^\top \Sigma_x a}{a^\top \Sigma_x a + \delta^2} \leq \varepsilon.
\end{align*}
Rearranging the above inequality results in $a^\top \Sigma_x a \leq \varepsilon \delta^2(1 - \varepsilon)^{-1}$ and, thus, the proof is complete. $\hfill \square$

\textit{Proof of Theorem 1.} The proof consists of four main derivations: \textbf{(i)} a quadratic expression for the value function, \textbf{(ii)} linear constraints guaranteeing hard input constraints, \textbf{(iii)} linear and quadratic inequalities for the softened versions of \eqref{eq:expectedvalueterm} and \eqref{eq:vartermCC}, respectively, and \textbf{(iv)} SOCP formulation of the soft-constrained version of the optimal control problem~\eqref{eq:P2}.

\textbf{(i)} Based on the exact penalty function method, the value function is adapted as $V_N(x,\mathbf{M},\mathbf{v}) + \rho \| \boldsymbol\epsilon \|_1$. Since by definition $\boldsymbol\epsilon \geq 0$, $\| \boldsymbol\epsilon \|_1 = \mathbf{1}^\top \boldsymbol\epsilon$ can be written as a linear constraint in the slack variables. Now, consider the value function \eqref{eq:valuefuncompact}. Substituting the system dynamics \eqref{eq:predmodelcompact} into \eqref{eq:valuefuncompact} and assigning $\tilde{x}_0 = x$ gives
\begin{align*}
\begin{array}{ll}
& V_N(x,\mathbf{u})  = \mathbf{E}[\| \mathbf{A}x + \mathbf{B}\mathbf{u} + \mathbf{D}\mathbf{G}\mathbf{w} \|^2_\mathbf{Q} + \| \mathbf{u} \|^2_\mathbf{R} ] \\
& = c(x) + \mathbf{E}[\|\mathbf{u} \|^2_{\mathbf{B}^\top\mathbf{Q}\mathbf{B} + \mathbf{R}} + 2\mathbf{u}^\top\mathbf{B}^\top\mathbf{Q}(\mathbf{A}x + \mathbf{D}\mathbf{G}\mathbf{w}) ],
\end{array}
\end{align*}
where $c(x) \triangleq \| x \|^2_{\mathbf{A}^\top \mathbf{Q}\mathbf{A}} + \mathbf{E}[\|\mathbf{w} \|^2_{\mathbf{G}^\top\mathbf{D}^\top\mathbf{Q}\mathbf{D}\mathbf{G}} ]$ since the states $x$ are known. Substituting the control policy \eqref{eq:affinedistcompactsat} into the above equation and rearranging yields 
\begin{align*}
\begin{array}{ll}
& V_N(x,\mathbf{\mathbf{M},\mathbf{v}}) =  c(x) + \mathbf{E}[ \|\mathbf{M}\boldsymbol\varphi(\mathbf{G}\mathbf{w}) + \mathbf{v} \|^2_{\mathbf{S}_1}]  \\
 &  + \mathbf{b}^\top\mathbf{E}[\mathbf{M}\boldsymbol\varphi(\mathbf{G}\mathbf{w}) + \mathbf{v}] +  \text{tr}\big(\mathbf{S}_2\mathbf{E}[(\mathbf{M}\boldsymbol\varphi(\mathbf{G}\mathbf{w}) + \mathbf{v})\mathbf{w}^\top] \big).
\end{array}
\end{align*}
By assumption, $\mathbf{E}[\boldsymbol\varphi(\mathbf{G}\mathbf{w})] = 0$ so that
\begin{align*}
\begin{array}{rl}
 V_N(x,\mathbf{\mathbf{M},\mathbf{v}})  = & c(x) + \| \mathbf{v} \|^2_{\mathbf{S}_1} + \mathbf{E}[ \| \mathbf{M}\boldsymbol\varphi(\mathbf{G}\mathbf{w}) \|^2_{\mathbf{S}_1} ] + \\
 & \mathbf{b}^\top\mathbf{v} + \text{tr}(\mathbf{S}_2\mathbf{M}\Omega_2) \\
= & c(x) + \mathbf{b}^\top\mathbf{v} + \\
& \| \mathbf{v} \|^2_{\mathbf{S}_1} + \text{tr}(\mathbf{M}^\top\mathbf{S}_1\mathbf{M}\Omega_1 + \mathbf{S}_2\mathbf{M}\Omega_2).
\end{array}
\end{align*} 
As $c(x)$ is independent of the decision variables $(\mathbf{M},\mathbf{v},\mathbf{Z},\boldsymbol\epsilon)$, it can be dropped from the value function. For matrices $A$, $B$, and $C$ of appropriate dimensions, the properties $\text{tr}(ABC)=\text{tr}(CAB)=\text{tr}(BCA)$, $\text{vec}(ABC) = (C^\top \otimes A)\text{vec}(B)$, and $\text{tr}(A^\top B) = \text{vec}(A)^\top \text{vec}(B)$ hold. Using these, the constraint $\mathbf{m} = \text{vec}(\mathbf{M})$, and the expressions for $\nu$ and $\Lambda$, the asserted value function can be obtained through standard algebraic manipulations. Note that the derived value function is a convex quadratic function since the matrices $\mathbf{S}_1\in\mathbb{S}_+^{n_u N}$ and $\Lambda \in \mathbb{S}_+^{n_u n_x N^2}$ are positive semidefinite.

\textbf{(ii)} From Lemma 1, it is known that the constraints $\mathbf{H}_u\mathbf{v} + \mathbf{Z}^\top \mathbf{k}_w \leq \mathbf{k}_u$, $\mathbf{Z}^\top \mathbf{H}_w = \mathbf{H}_u\mathbf{M}$, and $\mathbf{Z} \geq 0$ exactly represent $\mathbf{Pr}\big[\mathbf{H}_u\mathbf{u} \leq \mathbf{k}_u\big] = 1$ for the control policy \eqref{eq:affinedistcompactsat}. The latter constraints are linear in the decision variables $(\mathbf{M},\mathbf{v},\mathbf{Z})$ and, therefore, convex.

\textbf{(iii)} From the control policy \eqref{eq:affinedistcompactsat} and the assumption $\mathbf{E}[\boldsymbol\varphi(\mathbf{G}\mathbf{w})] = 0$, it can be derived that $\bar{\mathbf{u}} = \mathbf{v}$, $\Sigma_\mathbf{u} = \mathbf{M}\Omega_1\mathbf{M}^\top$, and $\sigma_{\mathbf{u},\mathbf{w}} = \mathbf{M}\Omega_2$. Substituting these definitions into \eqref{eq:xmeancompact} and \eqref{eq:Sigmaxcompact} yields $\bar{\mathbf{x}} = \mathbf{A}x_0 + \mathbf{B}\mathbf{v}$ and $\Sigma_\mathbf{x} = \mathbf{B}\mathbf{M}\Omega_1\mathbf{M}^\top\mathbf{B}^\top + \mathbf{B}\mathbf{M}\Omega_2\mathbf{G}^\top\mathbf{D}^\top + \mathbf{D}\mathbf{G}\Omega_2\mathbf{M}^\top\mathbf{B}^\top + \mathbf{D}\mathbf{G}\Sigma_\mathbf{w}\mathbf{G}^\top\mathbf{D}^\top$. Constraint \eqref{eq:softconsexp} is obtained by substituting the expression for $\bar{\mathbf{x}}$ into \eqref{eq:expectedvalueterm} and introducing the corresponding slack variables $\epsilon_{i,j}^m$. The latter constraint is linear in $(\mathbf{v},\epsilon_{i,j}^m)$ and, therefore, convex. Substituting the expression for $\Sigma_\mathbf{x}$ into \eqref{eq:vartermCC} results in $\mathbf{y}^\top_{i,j}\Omega_1\mathbf{y}_{i,j} + \mathbf{y}_{i,j}^\top\Omega_2\mathbf{G}^\top\mathbf{D}^\top\mathbf{H}^{[i]\top}_{x,j} + \mathbf{H}^{[i]}_{x,j}\mathbf{D}\mathbf{G}\Sigma_\mathbf{w}\mathbf{G}^\top\mathbf{D}^\top\mathbf{H}^{[i]\top}_{x,j} \leq \alpha_k\delta^2_k(1-\alpha_k)^{-1}$. It is evident that $\mathbf{y}_{i,j}^\top\Omega_2\mathbf{G}^\top\mathbf{D}^\top\mathbf{H}^{[i]\top}_{x,j} = (\mathbf{y}_{i,j}^\top\Omega_2\mathbf{G}^\top\mathbf{D}^\top\mathbf{H}^{[i]\top}_{x,j})^\top$ since this is a scalar constraint. Using the latter expression and definitions for $\mathbf{q}^\top_{i,j}$ and $c_{i,j}$, it can be stated that $\mathbf{y}^\top_{i,j}\Omega_1\mathbf{y}_{i,j} + \mathbf{q}_{i,j}^\top \mathbf{y}_{i,j} + c_{i,j} \leq 0$. Constraint \eqref{eq:softconsvar} is obtained by introducing the corresponding slack variables $\epsilon_{i,j}^v$. Note that \eqref{eq:softconsvar} is linear in $\epsilon_{i,j}^v$ and quadratic in $\mathbf{M}$ and, therefore, convex since $\Omega_1 \in \mathbb{S}_+^{n_xN}$.

\textbf{(iv)} A special case of SOCPs includes convex quadratically constrained quadratic programs. In the prequel, it was shown that all the constraints in~\eqref{e_P3} are linear or quadratic convex expressions. Problem \eqref{e_P3} is converted into the standard SOCP form by introducing a new decision variable $g$ that bounds the quadratic part of the value function, i.e., $\mathbf{v}^\top \mathbf{S}_1 \mathbf{v} + \mathbf{m}^\top \Lambda \mathbf{m} \leq g$. The value function in \eqref{e_P3} is then replaced by $\mathbf{b}^\top\mathbf{v} + \nu^\top\mathbf{m} + \mathbf{1}^\top \boldsymbol\epsilon + g$. Any convex quadratic constraint of the form $x^\top Q x + b^\top x + c \leq 0$ with $x \in \mathbb{R}^n$ and $Q \in \mathbb{S}^n_+$ can be formulated as an equivalent SOC constraint,
\begin{align*}
\left\| \begin{bmatrix}
(1+b^\top x + c)/2 \\
Q^{1/2} x
\end{bmatrix} \right\|_2 \leq (1-b^\top x -c)/2.
\end{align*}
Hence, Problem \eqref{e_P3} can be stated in the standard SOCP form by applying the above result to the new constraint $\mathbf{v}^\top \mathbf{S}_1 \mathbf{v} + \mathbf{m}^\top \Lambda \mathbf{m} \leq g$ and \eqref{eq:softconsvar}. $\hfill \square$

\textit{Proof of Lemma 4.} First, from the law of iterated expectations, $\mathbf{E}_{x_0}[V(x_t)] = \mathbf{E}_{x_0}[\mathbf{E}[V(x_t)|\{x_s\}_{s=1}^{t-1}]] = \mathbf{E}_{x_0}[\mathbf{E}_{x_{t-1}}[V(x_t)]]$. Next, it can be derived from the premises in Lemma 4 that $\mathbf{E}_{x_{t-1}}[V(x_t)] \leq \lambda V(x_{t-1})$ whenever $x_{t-1} \not\in \mathbb{D}$ and $\mathbf{E}_{x_{t-1}}[V(x_t)] \leq b$ whenever $x_{t-1}\in\mathbb{D}$. Adding these results together gives $\mathbf{E}_{x_{t-1}}[V(x_t)] \leq \lambda V(x_{t-1}) \mathbf{1}_{\mathbb{R}^{n_x}\backslash\mathbb{D}}(x_{t-1}) + b\mathbf{1}_{\mathbb{D}}(x_{t-1})$ for all $x_{t-1} \in \mathbb{R}^{n_x}$. Taking the expectation over $x_0$ and applying the first result yields $\mathbf{E}_{x_0}[V(x_t)] \leq \lambda \mathbf{E}_{x_0}[V(x_{t-1})] + b\mathbf{Pr}_{x_0}[x_{t-1} \in \mathbb{D}]$. Repeating these steps for $\{ \mathbf{E}_{x_0}[V(x_{s})] \}_{s=1}^{t-1}$ and recursively substituting the derived inequalities into the last inequality results in
\begin{align*}
\mathbf{E}_{x_0}[V(x_t)] & \leq \lambda^t V(x_0) + b\sum_{i=0}^{t-1} \lambda^{i} \mathbf{Pr}_{x_0}[x_{t-1-i} \in \mathbb{D}] \\
& \leq \lambda^t V(x_0) + b\sum_{i=0}^{t-1} \lambda^{i} \\
& \leq \lambda^t V(x_0) + b\sum_{i=0}^{\infty} \lambda^{i} \\
& \leq \lambda^t V(x_0) + b(1-\lambda)^{-1}.
\end{align*}
Note that $\mathbf{Pr}_{x_0}[x_{t-1-i} \in \mathbb{D}] \leq 1$ and the geometric series expression $1 + \lambda + \lambda^2 + \cdots = (1-\lambda)^{-1}$, which is convergent for $| \lambda | < 1$. The above result implies that $\sup_{t \in \mathbb{N}_0} \mathbf{E}_{x_0}[V(x_t)] \leq V(x_0) + b(1-\lambda)^{-1} < \infty$ is bounded and, thus, the proof is complete. $\hfill \square$

\textit{Proof of Theorem 2.} The proof follows directly from \cite{hokayem09}. Recall Holder's inequality $| a^\top b| \leq \|a\|_p \|b\|_q$ for all $p,q\in \mathbb{N}$ such that $p^{-1} + q^{-1} = 1$, $\|X Y \|_p \leq \|X\|_p \|Y\|_p$ for all $p \in \mathbb{N}$, and $\| a \|_\infty \leq \|a \|_1 \leq n\|a\|_\infty$ and $\|a\|_\infty \leq \| a \| \leq \sqrt{n}\|a\|_\infty$, where $a,b\in \mathbb{R}^n$, $X \in \mathbb{R}^{n\times m}$, and $Y \in \mathbb{R}^{m\times q}$. 

Since $A$ is assumed Schur stable (see Assumption 1), there exists a matrix $P \in \mathbb{S}_{++}^{n_x}$ such that $A^\top P A - P \leq - I_{n_x}$. Define the measurable function $V(x) = x^\top P x$. From \eqref{eq:CLrecursion}, it can be written that 
\begin{align*}
& \mathbf{E}_{x_t}[x_{t+1}^\top P x_{t+1}] = x_t^\top A^\top P A x_t + 2x_t^\top A^\top P B\kappa^s_N(x_t) \\
& ~~~~~~ + \kappa^s_N(x_t)^\top B^\top P B \kappa^s_N(x_t) + \text{tr}(G^\top P G \Sigma_w).
\end{align*}
As the hard input constraints \eqref{eq:inputcons} are assumed to be compact, there exists some set $\mathbb{U}_1 \triangleq \{ u \in \mathbb{R}^{n_u} \mid \| u \|_1 \leq U_b \}$ and $U_b \in \mathbb{R}_+$ such that $\mathbb{U} \subseteq \mathbb{U}_1$. Since $\kappa^s_N(x) \in \mathbb{U}$ must hold in $\mathbb{P}_N^s(x)$, $\| \kappa^s_N(x) \|_1 \leq U_b$ for all $x \in \mathbb{R}^{n_x}$. An upper bound can now be derived for the second and third terms in the above expression for $\mathbf{E}_{x_t}[x_{t+1}^\top P x_{t+1}]$, i.e., 
\begin{align*}
\begin{array}{ll}
2x_t^\top A^\top P B\kappa^s_N(x_t) &= 2(B^\top P A x_t)^\top \kappa^s_N(x_t) \\
& \leq 2 \|B^\top P A x_t\|_\infty \| \kappa^s_N(x_t) \|_1 \\
& \leq 2\|B^\top P A \|_\infty U_b \| x_t \|_\infty
\end{array}
\end{align*}
\begin{align*}
\begin{array}{ll}
\kappa^s_N(x_t)^\top B^\top P B \kappa^s_N(x_t) &= (B^\top P B \kappa^s_N(x_t))^\top\kappa^s_N(x_t) \\
& \leq \| B^\top P B \kappa^s_N(x_t) \|_\infty \| \kappa^s_N(x_t) \|_1 \\
& \leq \| B^\top P B \|_\infty \| \kappa^s_N(x_t) \|_1^2 \\
& \leq \| B^\top P B \|_\infty U_b^2.
\end{array}
\end{align*}
Define $c_1 = \|B^\top P A \|_\infty U_b$ and $c_2 = \| B^\top P B \|_\infty U_b^2 + \text{tr}(G^\top P G \Sigma_w)$ such that $\mathbf{E}_{x_t}[x_{t+1}^\top P x_{t+1}] \leq x_t^\top A^\top P A x_t + 2 c_1 \| x_t \|_\infty + c_2 \leq x_t^\top P x_t - \|x_t\|^2 + 2 c_1\| x_t \|_\infty + c_2$. By definition $c_1 \in \mathbb{R}_{++}$ and $c_2 \in \mathbb{R}_{++}$ such that $2 c_1\| x_t \|_\infty + c_2 \leq \theta \| x_t \|_\infty^2$ for some $\theta \in \mathbb{R}_{++}$. The latter expression can be rewritten as $-\theta(\| x_t \|_\infty -\frac{c_1}{\theta})^2 + c_2 + \frac{c_1^2}{\theta} \leq 0$ such that $\frac{1}{\theta}(c_1 + \sqrt{c_1^2 + c_2\theta}) \leq \| x_t \|_\infty$. Let $r = \frac{1}{\theta}(c_1 + \sqrt{c_1^2 + c_2\theta})$ and define $\mathbb{D} \triangleq \{ x\in \mathbb{R}^{n_x} \mid \| x \|_\infty \leq r \}$. For all $x_t \not\in \mathbb{D}$, the definition of $\mathbb{D}$ results in $2 c_1 \| x_t \|_\infty + c_2 \leq \theta \| x_t \|_\infty^2 \leq \theta \| x_t \|^2$ such that $\mathbf{E}_{x_t}[x_{t+1}^\top P x_{t+1}] \leq x_t^\top P x_t - (1-\theta) \|x_t\|^2$. As $P \in \mathbb{S}_{++}^{n_x}$, it is known that $x_t^\top P x_t \leq \lambda_{\text{max}}(P)\|x_t\|^2$ with $\lambda_{\text{max}}(P) > 0$ yielding $\mathbf{E}_{x_t}[x_{t+1}^\top P x_{t+1}] \leq (1 - \frac{1-\theta}{\lambda_{\text{max}}(P)})x_t^\top P x_t$ for all $\| x_t \|_\infty \not\in \mathbb{D}$. For the latter expression to represent a decay outside of $\mathbb{D}$, the constant $\theta$ must satisfy $1 - \lambda_{\text{max}}(P) < \theta < 1$. Since $\lambda_{\text{max}}(P) > 0$, there exists at least one $\theta$ that satisfies these inequalities.

Hence, the premises of Lemma 4 are satisfied with the following definitions: $V(x) \triangleq x^\top P x$, $\mathbb{D} \triangleq \{ x\in \mathbb{R}^{n_x} \mid \| x \|_\infty \leq r \}$, $b \triangleq \sup_{x_t \in \mathbb{D}}\mathbf{E}_{x_t}[x_{t+1}^\top P x_{t+1}]$, and $\lambda \triangleq (1 - \frac{1-\theta}{\lambda_{\text{max}}(P)})$ for any $\theta \in [1 - \lambda_{\text{max}}(P),1]$. It directly follows that $\{ \mathbf{E}_{x_0}[ x_{t}^\top P x_{t} ] \}_{t \in \mathbb{N}_0}$ is bounded for all $x_0 \in \mathbb{R}^{n_x}$. Since $\lambda_{\text{min}}(P)\|x_t\|^2 \leq x_t^\top P x_t$, the assertion holds and the proof is complete. $\hfill \square$

\section*{Appendix B}
\label{sec:app_B}

This appendix contains the description of the state-space model of the continuous acetone-butanol-ethanol fermentation process used in Section \ref{sec:case_study}. The state vector is defined as
\begin{align*}
\begin{array}{ll}
x = &[C_{AC} \; C_{A} \; C_{En} \; C_{AaC} \; C_{Aa}  \; C_{BC} \; C_{B} \\
 & C_{An} \; C_{Bn} \; C_{Ad} \; C_{Cf} \; C_{Ah}]^\top,
\end{array}
\end{align*}
where $C$ denotes concentration (mM) of AC = Acetyl-CoA, A = Acetate, En = Ethanol, AaC = Acetoacetate-CoA, Aa = Acetoacetate, BC = Butyryl-CoA, B = Butyrate, An = Acetone, Bn = Butanol, Ad = adc, Cf = ctfA/B, and Ah = adhE (see \cite{haus11}). The input vector is defined as
\begin{align*}
u = [D \; G_0]^\top,
\end{align*}
where $D$ is the dilution rate ($\text{hr}^{-1}$); and $G_0$ is the inlet glucose concentration (mM). The system matrices are given by 
{\tiny
\begin{align*}
& A = \\
& \left [ \arraycolsep=1pt\begin{array}{cccccccccccc}
0.888& 0.017& 0& 0.143& 0& 0& -0.002& 0& 0& 0& 0.017& -0.002\\
-0.001& 0.969& 0& -0.151& 0& 0& 0.002& 0& 0&0& -0.018& -0.002\\
0.083& 0.001& 0.988& 0.007& 0& 0& 0& 0& 0& 0& 0.001& 0.002\\
0.015&  -0.017& 0&    0.716& 0& 0&  -0.018& 0& 0& 0&   -0.032& 0\\
0& 0& 0& 0& 0& 0& 0& 0& 0& 0& 0& 0\\
0& 0& 0& 0& 0& 0& 0& 0& 0& 0& 0& 0\\
-0.001& 0.001& 0&  -0.122& 0&  0&  0.968&  0& 0& 0& -0.014& 0\\
0.003&  0.017& 0&  0.273&   0.988&  0& 0.018& 0.988&  0& 0& 0.0321& 0\\
0& 0& 0& 0& 0& 0& 0& 0& 0& 0.988& 0& 0\\
0& 0& 0& 0& 0& 0& 0& 0& 0& 0& 0.988& 0\\
0& 0& 0& 0& 0& 0& 0& 0& 0& 0& 0& 0.988 \end{array} \right ], 
\end{align*}
\begin{align*}
& B^\top = -\left [ \arraycolsep=1.4pt\begin{array}{cccccccccccc}
0.220& 2.39& 1.23& 0.157& 0& 0& 1.86& 7.33& 8.63& 0.242& 2.34& 6.29\\
0& 0& 0& 0& 0& 0& 0& 0& 0& 0& 0& 0 \end{array} \right ], 
\end{align*}
}
and $G = \text{diag}(1,1,1,1,1,1,1,1,1,1,1,1)$.

\bibliographystyle{elsarticle-num}
\bibliography{Literature_list}





\end{document}